\documentclass{article}
\usepackage{amssymb}
\usepackage{graphicx}
\usepackage{amsmath}

\newtheorem{proposition}{Proposition}
\newtheorem{lemma}{Lemma}

\def\NN{{\mathbb{N}}}

\def\CC{{\mathcal{C}}}

\def\II{{\mathcal{I}}}
\def\JJ{{\mathcal{J}}}
\def\TT{{\mathcal{T}}}

\def\SS{{\mathcal{S}}}
\def\WW{{\mathcal{W}}}
\def\KK{{\mathcal{K}}}
\def\MM{{\mathcal{M}}}

\title{Dual heuristics and new dual bounds to schedule the maintenances of nuclear power plants}
\author{Nicolas Dupin, El-Ghazali Talbi\\
Univ. Lille, UMR 9189 - CRIStAL \\
              {nicolas.dupin.2006@polytechnique.org}}

\date{
}

\begin{document}


\maketitle

\noindent{\textbf{Abstract}:
The EURO/ROADEF 2010  Challenge aimed to schedule the maintenance and refueling operations of French nuclear power plants,
 ranking the  approaches in competition for the quality of primal solutions. 
This paper  justifies the high quality of the best solutions computing  dual bounds with dual heuristics. 
 A first step designs several Mixed Integer  Programming (MIP) 
relaxations with different compromises between computation time and  quality of dual bounds. 
To deal with smaller  MIPs,  we prove how reductions
 in the number of time steps and  scenarios can guarantee dual bounds for the whole problem of the Challenge.
 Several sets of dual bounds are computable,  improving significantly the former best dual bounds of the literature.
 Intermediate results  allow also a better understanding of the problem 
 and offer  perspectives to improve some approaches of the Challenge.\\
\textbf{Keywords}: { Mixed Integer Programming ; Stochastic Programming ; Dual bounds 
; EURO/ROADEF 2010 Challenge ; Maintenance scheduling}
}

\section{Introduction}

The  ROADEF/EURO 2010 Challenge  was  specified  by  the French utility company (EDF) in \cite{roadef}
to  to schedule the maintenance and refueling of nuclear power plants. 
This is not a pure scheduling problem: 
the maintenance scheduling is optimized regarding the expected production costs to fulfill the power demands among the available power plants,
taking also into account the technical constraints of the power generation.
The optimization problem was formulated using 2-stage stochastic programming for the challenge.
Power demands, production capacities and costs are stochastic,  discrete scenarios model this uncertainty. 

The  contributions of the challengers gave rise to a special issue of Journal of Scheduling \cite{ozcan2013}.
A major difficulty of this Challenge is that it is a
complex industrial optimization problem with varied constraints and types of variables.
The best results were mainly obtained in \cite{gardi} with an aggressive local search approach.
Mixed Integer  Programming (MIP) formulations appeared naturally with many linear constraints in the specification \cite{roadef}.
Several approaches in competition were based on exact methods like (\cite{jost,lusby,Roz12}), this required heuristic reductions,
constraint relaxations and repairing procedure to compute final solutions.
The large size of the instances was a bottleneck for the exact methods.

It was an open question after the Challenge to have dual bounds for this large scale problem.
The only  dual bounds  published were \cite{Bra13} with optimal computations of very relaxed problems. 
 This paper investigates how to compute efficiently dual bounds for the problem of the Challenge
using dual heuristics for MIP as defined in \cite{li2015dual}: 
in Linear Programming (LP) based Branch-and-Bound (B\&B) algorithms, 
many heuristics have been developed to improve the efficiency of B\&B  search of primal solutions, 
but on the dual side, B\&B relied solely on cutting planes to improve dual bounds. 
Dual heuristics improve dual bounds with other relaxation algorithms. 
New dual bounds for the ROADEF Challenge are provided in this paper with restricted MIP  computations.
A crucial point will be to prove that the restricted MIP computations
induce dual bounds for the original problem.

This paper is organized as follows. 
Section 2 gives an overview of the problem constraints, while section 3 presents 
the state-of-the art of the solving approaches. 
Section 4 proposes a MIP model for the problem, relaxing only two sets of constraints (CT6 and CT12).
Section 5 proposes a light formulation of constraints CT6 to provide better dual bounds.
Section  6 proposes some parametric relaxations. 
Section  7 analyzes how the aggregation  of production time steps or scenarios
allow to compute dual bounds for the original problem. 
The computational analyses and results are reported in section 8.
The conclusions and perspectives of this work  are drawn in section 9.

\begin{figure}[ht]
      \centering
            \includegraphics[angle=0, width=.87\linewidth]{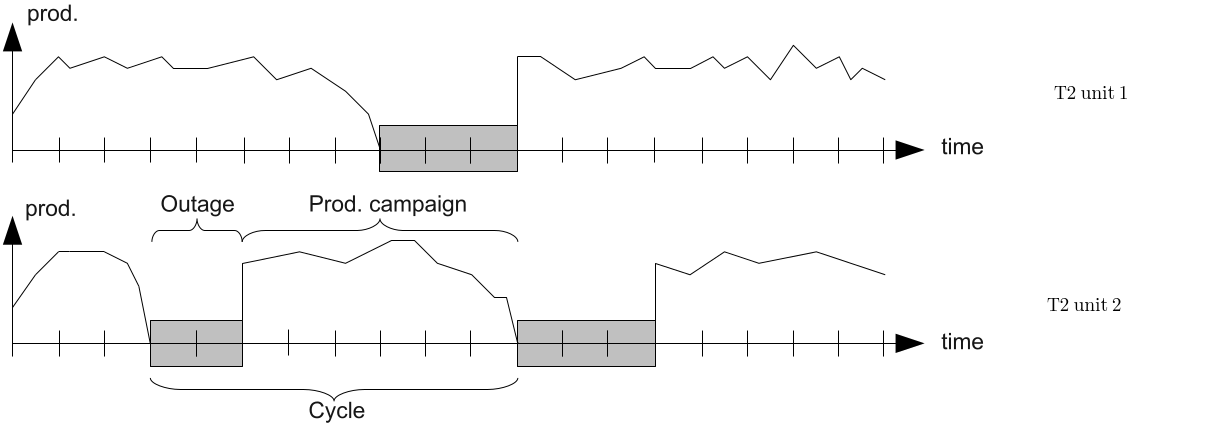}
      \caption{Illustration of definitions: cycles, production campaigns and outages}\label{defIllustr}  
\end{figure}


 \begin{table}
\noindent{
\caption{Notations of the set and indexes}\label{setROADEF}
\begin{tabular}{ll}
$t\in \TT = [\![1, T]\!]$ & Production time steps, index $t$ corresponds to period $[t,t+1]$\\  
$w \in \WW = [\![1, W]\!]$ & Weekly time steps to place outage dates.\\
$j \in \JJ = [\![1, J]\!]$ & Flexible (Type 1, T1) power plants.\\
$i \in \II = [\![1, I]\!]$ & Nuclear power plants (Type 2, T2).\\
$k \in \KK = [\![0, K]\!]$ & Cycles related to T2 units, $k=0$  for initial conditions.\\ 
$s\in \SS = [\![1, S]\!]$ & Stochastic scenarios for demands, production costs and capacities.\\
$m \in \MM_{i,k}= [\![1,\mathbf{Np}_{i,k}]\!]$ & Points of the  CT6 profile of T2 plant $i$ ending the cycle $(i,k)$.
\end{tabular}
}
 
 \vskip 0.3cm
 
 \caption{Definitions and notations for the imput parameters}\label{notationsROADEF}

 \textbf{Temporal notations }

\noindent{
\begin{tabular}{ll}
$\mathbf{Dem}^{t,s}$ & Power demands at time step $t$ for the scenario $s$.\\
$\mathbf{D}_{t}$ & Duration between time steps of indexes  $t$ and $t+1$.\\
$\mathbf{F}_{t}$ & Conversion factor between power and fuel in time step $t$.\\
  $t_w$ & Index of the first time step of week $w$.\\
 $w_t$ & Week of the production time step  $t$.
\end{tabular}
}

\vskip 0.3cm

\textbf{Notations for T1 units $j$}

\noindent{
\begin{tabular}{ll}
$\mathbf{C}^s_{jt}$ & Production Costs  proportional to the generated power at $t$ for scenario $s$.\\
$\mathbf{\underline{P}}^s_{jt}$ & Minimal power to generate at time step $t$ for scenario $s$. \\
$\mathbf{\overline{P}}^s_{jt}$ & Maximal generated  power at time step $t$ for scenario $s$. 
\end{tabular}
}

\vskip 0.3cm

\textbf{Notations for T2 units $i$}

\begin{tabular}{ll}
$\mathbf{\overline{A}}_{i,k}$ & Maximal fuel level remaining in cycle $k$ to process outage  $k+1$.\\
$\mathbf{Bo}_{i,k}$ & Fuel level "Bore O" of cycle $k$\\
$\mathbf{C}_i$  & Proportional cost to the final remaining fuel levels at $W$.\\
$\mathbf{c}_{i,k,m}$ & Loss coefficient of generated power at mode $m$ and cycle $k$.\\
$\mathbf{Da}_{i,k}$ & Outage duration for maintenance and refueling at cycle $k$.\\
$\mathbf{{To}}_{i,k}$ & First possible outage week for cycle $k$ of T2 plant $i$.\\
$\mathbf{{Ta}}_{i,k}$ & Last possible beginning  week for outage $k$ of T2 plant $i$.\\
$\mathbf{f}_{i,k,m}$ & Fuel levels for the stretch decreasing profile at mode $m$ and cycle $k$.\\
$\mathbf{Mmax}_{i,k}$ & Maximal modulation  for cycle $k$.\\
$\mathbf{\overline{P}}_i^t$ & Maximal generated  power at time step $t$ .\\
$\mathbf{Q}_{i,k}$ & Proportion of fuel that can be kept during reload in cycle $k$ at plant $i$\\
$\mathbf{\underline{R}}_{ik}$ & Minimal refueling at outage $k$.\\
$\mathbf{\overline{R}}_{ik}$ & Maximal refueling at outage $k$.\\
$\mathbf{\overline{S}}_{i,k}$ & Maximal fuel level  of T2 plant $i$ at   production  cycle $k$.\\
$ \mathbf{Xi}_{i}$ & Initial fuel  stock of T2 plant $i$.
\end{tabular}

 \end{table}
 
\section{Problem statement}

This section presents the model specified for the EURO/ROADEF 2010  Challenge.
 Mathematical notations are gathered in Tables \ref{setROADEF} and \ref{notationsROADEF}.

\subsection{Physical assets and time step discretization}


\paragraph{\textbf{Production assets and maintenances}}
To generate electricity, two kinds of power plants are modeled for the EURO/ROADEF 2010  Challenge. 
On one hand, Type-1 (shortly T1, denoted $j \in \JJ$) plants model coal, fuel oil, or gas power facilities. 
T1 power plants can be supplied in fuel continuously without inducing any  offline periods. 
On the other hand, Type-2 (shortly T2, denoted $i \in \II$) power plants have to be shut down for refueling
and maintenance regularly. T2 units correspond to nuclear power plants. 

 The production planning for a T2 unit is organized in a succession of cycles, an offline period (called \emph{outage}) 
 followed by an online period (called \emph{production campaign}), as illustrated in Figure \ref{defIllustr}. 
Cycles are indexed with $k \in \KK = [\![0, K]\!]$, each  T2 unit $i$ having possibly $K$ maintenances scheduled in the time horizon ($K$ is common for all T2 units in the challenge).
The number of maintenances planned in the time horizon may be inferior to $K$.
Cycle $k=0$ denotes initial conditions, the current cycle at $w=0$. 
For  units on maintenance at $w=0$, cycle $k=0$ considers the remaining duration of the outage.
For  units on production at $w=0$, notations and constraints are extended considering  a fictive cycle $k=0$, with a null duration of the corresponding outage.
For both cases, the  refueling of the cycle $k=0$ is equal to the initial fuel stock,
and the beginning week of this cycle is fixed.

\paragraph{\textbf{From maintenance scheduling to production planning}}
The time horizon is discretized with two kind of homogeneous time steps.
On one hand, outage decisions are discretized weekly and indexed with $w \in \WW=[\![1;W]\!]$.
On the other hand, production time steps for T1 and T2 units are discretized with  
$t \in \TT$, production times steps  from 8h to 24h.
This fine discretization is required to consider fluctuating demands in hourly periods, as emphasized in \cite{ucp}.

Maintenance scheduling and production planning are firstly coupled with the null production phases implied by outages.
Secondly, the production constraints of T2 power plants must be feasible once the outage decisions are made.
Another heterogeneity of maintenance and production  decisions is modeled: production decisions are considered for
a discretized set of  scenarios $s\in \SS$ to model uncertainty on power  demands and production capacities/costs.
The outage and refueling decisions are deterministic, and must guarantee the feasibility of the production planning for each scenario.

\subsection{Decision variables and objective function}

The outage  and production decisions are optimized conjointly.
The variables for the outage dates and the refueling levels are related to T2 plants,
whereas and the  production variables concern  T1 and T2 plants for all the time periods and for all the scenarios.

The objective function aims  to minimize the production costs.
T1 Production costs are proportional to the generated power. 
T2 production costs are proportional to the refueling quantities over the time horizon,
 minus the expected stock values at the end of the period to avoid end-of-side effects.
The global cost  is the average value of the production costs 
considering all the scenarios, following the 2-stage Stochastic Programming framework.
 
%


\subsection{Description of the constraints}
There are 21 sets of constraints in the Challenge,
 numbered from CT1 to CT21 in the  specification \cite{roadef}.
These constraints can be classified into two categories:
CT1 to CT12 concern the production constraints in the operational level whereas CT13 to CT21
are scheduling constraints in the strategic level.

\begin{figure}[ht]
      \centering
      \includegraphics[angle=0, width=.79\linewidth]{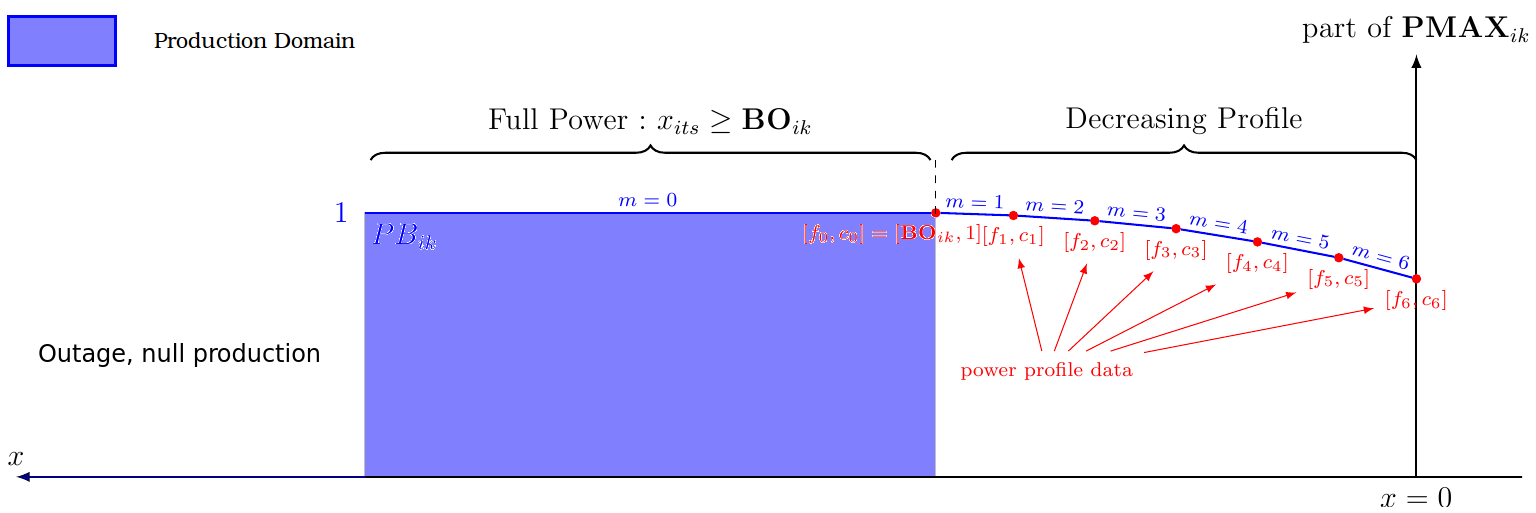}
      \caption{Production domain for nuclear units: null production during outages, flexible production
      in production campaign while the stock is upper $\mathbf{Bo}_{i,k}$, and then an imposed decreasing profile}\label{prodRoadef} 
\end{figure}

\paragraph{\textbf{Production constraints}}
Table \ref{contraintesROADEF} lists the production constraints.
CT1 couples the production of the T1 and T2 plants with power demands to fulfill exactly.
CT2 and CT3 induce a continuous domain for T1 productions.
CT4 to CT6 describe the production possibilities for T2 units depending on the fuel stock with a null production 
for the maintenance periods. 
As  illustrated in Figure \ref{prodRoadef}, 
the T2 production is deterministic  when the fuel stock is lower than a given threshold following a decreasing
piecewise affine function. 
CT7 to CT11 involve fuel constraints for T2 units: bounds on fuel stocks, on fuel refueling,
and relations between remaining fuel stocks and production/refueling operations.
CT12 imposes furthermore that T2 plants have limited capacities for modulations (i.e. for a non maximal production).

\begin{table}
\caption{Constraints coupling productions, fuel stock levels and outage decisions }\label{contraintesROADEF}

\begin{tabular}{l}
\hline
 \begin{minipage}{0.99 \linewidth}
 \vskip 0.2cm
 \textbf{CT1, demand covering}:  for all production time step $t\in \TT$ and all scenario $s \in \SS$, 
 the total production of T1 and T2 power plants must equalize the demands $\mathbf{Dem}^s_{t}$. 
   \vskip 0.3cm
 \end{minipage}\\
 \hline
  \begin{minipage}{0.99\linewidth}
 \vskip 0.2cm
 \textbf{CT2, T1 production bounds} for all production time step $t\in \TT$ and all scenario $s \in \SS$,
 the production domain of T1 plant $j\in \JJ$ describes exactly the continuous domain   $[\mathbf{\underline{P}}^s_{jt}, \mathbf{\overline{P}}^s_{jt}]$ 
 \vskip 0.3cm
 \end{minipage}\\
 \hline
 \begin{minipage}{0.99\linewidth}
 \vskip 0.2cm
  \textbf{CT3, null productions during outages}:
  for all scenario $s \in \SS$, the productions of  offline T2 power plants are null.
  \vskip 0.3cm
 \end{minipage}\\
\hline
 \begin{minipage}{0.99\linewidth}
 \vskip 0.2cm
  \textbf{CT4, CT5,  T2 production bounds}
  For all period $t\in \TT$ and all scenario $s \in \SS$,
 the production domain of T1 plant $j\in \JJ$ describes exactly the continuous domain   $[0,\mathbf{\overline{P}}_{it}]$ when its fuel level is superior to $\mathbf{Bo}_{i,k}$.
     \vskip 0.3cm
 \end{minipage}\\
\hline
 \begin{minipage}{0.99\linewidth}
 \vskip 0.2cm
  \textbf{CT6, "stretch" constraints}  (relaxed fully or partially in this study)
  During every scenario $s \in \SS$  and every time step $t \in \TT $ of the production campaign of cycle $k \in K_i$,
if the current fuel stock of plant $i \in \II$ is inferior to the level $\mathbf{Bo}_{i,k}$ , the production of $i$ is deterministic,
following a decreasing profile,  piecewise linear function of the stock level, as  in Figure \ref{prodRoadef}.  
   \vskip 0.3cm
 \end{minipage}\\
\hline

  \begin{minipage}{0.99\linewidth}
 \vskip 0.2cm
 \textbf{CT7, refueling quantities}:  
 The refueling possibilities for outage $k$ of T2 plant $i$ describes exactly the continuous domain $[\mathbf{\underline{R}}_{i,k},\mathbf{\overline{R}}_{i,k}]$
 \vskip 0.3cm
 \end{minipage}  
  \\
\hline
    \begin{minipage}{0.99\linewidth}
 \vskip 0.2cm
  \textbf{CT8, initial fuel stock}: The initial fuel stock for T2 unit $i$ is $\mathbf{Xi_i}$, known and common for all scenario $s$.
 \vskip 0.3cm
 \end{minipage}\\
\hline

      \begin{minipage}{0.99\linewidth}
 \vskip 0.2cm
  \textbf{CT9, fuel stock variation during a production campaign}: The
  fuel stock variation during a production campaign of a cycle $(i,k)$ between $t$ and $t+1$, 
  is  proportional to the power produced by $i$ at time step $t$, with a proportional factor $- \mathbf{F}_t$.
   \vskip 0.3cm
 \end{minipage}\\
\hline

\begin{minipage}{0.99\linewidth}
 \vskip 0.2cm
  \textbf{CT10, fuel stock variation during an outage}
  During an outage, the fuel stock variation is the sum of the decisional refueling bounded with CT7 and
  a certain amount of unspent fuel, calculated with
  a proportional loss with factor $\mathbf{Q}_{i,k}<1$ to the residual fuel before refueling.
   \vskip 0.3cm
 \end{minipage}\\
\hline 
    
   \begin{minipage}{0.99\linewidth}
 \vskip 0.2cm
 \textbf{CT11, bounds on fuel stock} The fuel level is in $[0,\mathbf{S}_{i,k}]$ for cycle $k$ of T2 unit $i$.
   The fuel level  must be lower than  $\mathbf{A}_{i,k+1}$ to process outage $k+1$. 
 \vskip 0.3cm
 \end{minipage}\\
 \hline
  \begin{minipage}{0.99\linewidth}
 \vskip 0.2cm
 \textbf{CT12, Modulation constraints} (relaxed fully in this study)
Production cycle of T2 unit $i$ is mostly the maximal powers $\mathbf{\overline{P}}_{it}$ when the fuel level is superior to $\mathbf{Bo}_{i,k}$ for technical constraints.
The deviations are limited to quantities $\mathbf{\overline{M}}_{ik}$ for all cycle $k$.
\vskip 0.3cm
 \end{minipage}\\
 \hline
 \end{tabular}
 
 \end{table}

 \paragraph{\textbf{Scheduling constraints}}
Table \ref{contraintesROADEFordo} lists the scheduling constraints.
The maintenances follow the order of set $k\in \KK$ without skipping maintenances: 
if an outage $k+1$ is processed in the time horizon, it must
follows the production cycle $k$.
CT13 impose time windows to begin some maintenances, which can  impose some maintenance dates.
CT14 to CT18 constrain minimum spacing or maximum overlapping among a subset of outages.
CT19 restrict simultaneous maintenances because of some limitations in resource usages (specialized tools or teams).
CT20 impose a maximum number of simultaneous outages 
while CT21 limits the maximum cumulated offline T2 power due to maintenances operations.

 \begin{table}
\caption{Scheduling constraints coupling the dates of outages}\label{contraintesROADEFordo}

\begin{tabular}{l}
\hline
  \begin{minipage}{0.99\linewidth}
 \vskip 0.2cm
  \textbf{CT13, time windows for maintenances}   constraints for the beginning dates of outages 
  \vskip 0.1cm
 
 \textbf{Implicit constraints}:  
 The maintenances follow the order of set $k\in \KK$ without skipping maintenances: if  outage $k+1$ is processed, it must
 follow the cycle $k$.
 \vskip 0.3cm

 \end{minipage}
  \\
\hline
  \begin{minipage}{0.99\linewidth}
 \vskip 0.2cm
  \textbf{CT14, minimal spacing/ maximal overlapping constraints}: 
     For all constraint $c\in \CC_{14}$,  a subset of outages, $A_c^{14}$  have to be spaced by at least $S_c^{14}$ weeks:
     for  $S_c>0$, it is a minimal spacing from the beginning of a previous outage to the beginning of a next outage,
     for  $S_c\leqslant 0$, it is a maximal number of weeks where two outages can overlap.
 \vskip 0.3cm
 \end{minipage}\\
\hline
  \begin{minipage}{0.99\linewidth}
 \vskip 0.2cm
   \textbf{CT15,  spacing/overlapping constraints during a specific period}: 
CT15 are similar to CT14 with sets $\CC_{15},A_c^{15}$ and parameters $S_c^{15}$,
applying only on a time interval $[U_c ; V_c ]$.
  \vskip 0.3cm
 \end{minipage}\\
\hline
  \begin{minipage}{0.99\linewidth}
 \vskip 0.2cm
  \textbf{CT16, Minimum spacing constraint between decoupling dates}: 
    For all constraint $c\in \CC_{16}$, a   subset of outages $A_c^{16}$
   have to be spaced by at least $S_c^{16}>0$ weeks from the beginning of the previous outage to the beginning of the next outage. 
  \vskip 0.3cm
 \end{minipage}\\
\hline
  \begin{minipage}{0.99\linewidth}
 \vskip 0.2cm
  \textbf{CT17, Minimum spacing constraint between dates of coupling}: 
  For all constraint $c\in \CC_{17}$, a   subset of outages $A_c^{17}$
   have to be spaced by at least $S_c^{17}>0$ weeks from the end of the previous outage to the end of the next outage. 
\vskip 0.3cm
 \end{minipage}\\
\hline

  \begin{minipage}{0.99\linewidth}
 \vskip 0.2cm
  \textbf{CT18, minimum spacing constraint between coupling and decoupling dates} : 
  For all constraint $c\in \CC_{18}$, a   subset of outages $A_c^{18}$
 have to be spaced by at least $S_c^{18}>0$ weeks from the end of the previous outage to the beginning of the next outage. 
\vskip 0.3cm
 \end{minipage}\\
\hline

  \begin{minipage}{0.99\linewidth}
 \vskip 0.2cm
  \textbf{CT19, resource constraints for maintenances}:
For all constraint $c\in \CC_{19}$, 
a subset of outages $(i_0,k_0) \in A_c^{19}$ shares a common resource (which can be a single specific tool or maintenance team)
with forbid simultaneous usage of this resource.
$U_{i_0,k_0,m_0}$ and $V_{i_0,k_0,m_0}$
 indicate respectively  the start and the length of the resource usage. 
\vskip 0.3cm
\end{minipage}\\
 \hline
 
 \begin{minipage}{0.99\linewidth}
 \vskip 0.2cm
 \textbf{CT20,  maximal number of simultaneous outages} 
For all constraint $c\in \CC_{20}$, 
at most $N_c (w)$ outages of of a subset of T2 plants $\II_c$ can overlap during the week $w\in \WW$.

 \vskip 0.3cm
 \end{minipage}\\
 \hline

 \begin{minipage}{0.99\linewidth}
 \vskip 0.2cm
 \textbf{ CT21, maximal power off-line} 
For all constraint $c\in \CC_{21}$, a time period $[U_c ; V_c ]$ is associated where the offline power  of T2 plants $\II_c\subset \II$
 due to maintenances must be lower than $I_c^{max}$. 
 \vskip 0.3cm
 \end{minipage}\\
 \hline
 \end{tabular}
 \end{table}

%

\section{Related work}

This section  describes the solving approaches for the EURO/ROADEF Challenge,
focusing on the dual bounds possibilities. 
We refer to \cite{froger2015maintenance} for a general survey on maintenance scheduling in the electricity industry, 
and to \cite{dupin2015modelisation} for a specialized review on the approaches for the 2010 EURO/ROADEF Challenge.



\subsection{\textbf{General facts}}
 We distinguish three  types of approaches for the  ROADEF challenge. 
The first category gathers \textit{MIP-based exact approaches}  using MIP models to tackle a simplified and reduced MIP problem
before a reparation procedure to build solutions for the original problem, as in \cite{jost,lusby,Roz12}.
 The second category gathers \textit{frontal local search approaches}. 
In this category, we find the two best approaches of the Challenge in terms of solution quality, \cite{gardi}
for an  aggressive Local search with the methodology of \cite{benoist2011localsolver},
\cite{beImproved} having similar results with a Simulated Annealing approach.
 The third category,  as in \cite{Ang12,brandt,Gav13,Godskesen}, gathers \textit{heuristic decomposition approaches}
solving iteratively   two types of subproblems.
The high-level  problem schedule outage dates and refueling decisions with constraints CT7-CT11 and CT13-CT21 
can be solved with MIP or Constraint Programming approaches. 
Low-level subproblems optimize independently for all scenario the production plans fulfilling the constraints CT1-CT6 and CT12
with fixed fuel levels and maintenance dates.



%

\subsection{\textbf{Reductions by preprocessing}}
Facing the large size instances of the Challenge, many approaches of the Challenge used preprocessing 
to reduce the size of the problem.
In that goal, exact processing can fix implied or necessarily optimal decisions,
which is a valid preprocessing to compute dual bounds.
Heuristic preprocessing strategies were also commonly used. 

\paragraph{\textbf{Exact preprocessing of time windows}}
Exact processing can be done tightening time window constraints CT13. 
To tighten exactly CT13 constraints, 
minimal durations  $\underline{D}_{i,k}$ of production campaigns $(i,k)$ are implied by
the initial stocks defined in CT8, maximal T2 productions CT5, max stock before refueling with CT11.
 Using this a preprocessing, some outages can be removed when their earliest completion time exceed the considered time horizon.


\paragraph{\textbf{Heuristic preprocessing of time windows}}
In the challenge specification, there is no maximal length on a production cycle: 
decreasing profile phases can be followed by a nil production phase without maximal duration in this state.
However, such situations are not economically profitable as the nuclear production cost is lower than the thermal production cost. 
This encouraged \cite{lusby} to avoid such situations to tighten the latest dates to begin outages $\mathbf{{Ta}}_{i,k}$.


\paragraph{\textbf{Reductions by aggregations}}
One of the main difficulties of the challenge is the big size of the instances.
It was a natural and common idea to reduce the sizes of the instances   with  aggregations. 
On one hand,   the power  production in the time steps $t\in \TT$ can be aggregated  to their weekly average value. 
On the other hand,  the stochastic scenarios were often  aggregated into one deterministic scenario considering the average 
values of T1 power bounds and production costs.
In order to reduce the number of scenarios, \cite{Gav13} focused on 3 scenarios including the average scenario and the scenarios with the minimal and maximal accumulated demands. 
Their tests conclude that it is a good compromise of the scenarios reduction which does not degrade the quality of the solutions.
None of the previous work studied the impact of such aggregations in terms of dual bounds.

\subsection{\textbf{MIP formulations}}
Several approaches of the Challenge were based on MIP formulations.
Many constraints are indeed formulated with linear formulation in the specification \cite{roadef}.
Generally, MIP models  tackled  simplified  MIP problems,
in a hierarchic heuristic repairing the previously relaxed constraints.
There were three types of MIP solving: straightforward B\&B searches using a compact formulation,
a Bender's decomposition and an extended formulation for a  column generation approach.

\paragraph{\textbf{Compact formulations}} 
Several MIP based approaches noticed that
scheduling constraints from CT14 to CT21 are modeled efficiently with MIP formulations
using time indexed formulation of constraints.
Jost and Savourey (2013)  used only these constraints in a two stage heuristic, computing production and stock levels
for solutions of the high level scheduling MIP with constraints CT13 to CT21 \cite{jost}.
Such formulation cannot provide dual bounds for the Challenge ROADEF because it does not model the production costs.
Solving MIP problems with scheduling constraints CT13 to CT21 and only binaries for outage dates in
\cite{jost} induce fast resolutions for all the instances of the Challenge.
%

We note that \cite{Jonc10} furnished as preliminary work  a MIP model for all the constraints. 
It is assumed that this formulation won't be able to tackle the large size of the instances.
This compact MIP formulation considered binaries for the outage decisions, but also binaries to follow exactly the decreasing profile
imposed with CT6 constraints.
Such variables were defined for all cycle $i,k$, for all $m \in \MM_{i,k}$, and for all $(s,t)\in \SS \times \TT$.

\paragraph{\textbf{Bender's decomposition approach}} 
Lusby et al. (2013)  provided the only  approach based on exact methods  which did not aggregate the stochastic scenarios \cite{lusby}.
Their  model relaxes fully the constraints CT6 and CT12, for a MIP formulation with binaries only for the outage decisions.
They aggregated time steps $t\in\TT$ into weeks for size reasons.
Their approach tackles the 2-stage stochastic programming structure  using a  Bender's decomposition: 
the master problem concerns the dates of outages and the refueling quantities, 
whereas  independent sub-problems are defined for each stochastic scenarios 
with continuous variables for productions and fuel levels.
Because of the limitations in  memory usage and time limit imposed for the Challenge, 
the Bender's decomposition algorithm is not deployed entirely. 
The  heuristic of \cite{lusby} computes first the LP relaxation exactly with the Bender's decomposition algorithm.
Then, a cut\&branch approach repairs integrity, branching on binary variables without adding new Bender's cuts.
The resulting heuristic approach was efficient for the small dataset A, difficulties and inefficiencies occur for the
 real sizes instances of the datasets B and X.
 If the aggregation of time steps lead to dual bounds for the original problem (not proven in \cite{lusby}, this paper will give some answers),
  dual bounds can be computed after the LP relaxation but also with the remaining dual bounds after the cut\&branch phase.
 
\paragraph{\textbf{Column generation approach}}
Rozenkopf et al. (2013)considered an exact formulation of CT6 constraints,
in a Column Generation (CG) approach dualizing coupling constraints among units \cite{Roz12}. 
 Such decomposition  requires to consider a single scenario, they also consider the average scenario.
Production time steps were also aggregated weekly.
These two simplifications are not prohibitive to compute dual bounds (this article will prove it).
A third simplification is prohibitive to compute dual bounds: the production domains are discretized 
to solve CG subproblems by dynamic programming.
This is a heuristic reduction of the feasible domain, 
dual bounds are computed on a heuristically reduced problem.
CG approach is deployed to compute a LP relaxation relaxing scheduling constraints CT14-CT21.
The further CG heuristic incorporate these scheduling constraints in the integer resolution with the columns generated, 
without adding more columns thereafter and without further branching heuristics.
Time and memory limitations explain such heuristic reduction.

\subsection{\textbf{Dual bounds}}

It was an open question after the Challenge to have dual bounds for this large scale problem.
As mentioned previously, none of the exact methods derived dual bounds for the  ROADEF problem
for the real size instances.
Semi-definite programming relaxations were also investigated in \cite{Gor12},
but the size of the instances is still a bottleneck.
Brandt et al. published the only dual bounds for the whole problem of the EURO/ROADEF 2010 Challenge \cite{Bra13}.
Two methods were investigated, computing optimal solutions of highly relaxed problems. 
Their first  method relaxes power profile constraint (CT6), as well as all fuel
level tracking (CT7 to CT12) and outage scheduling constraints (CT13 to CT21).
The remaining computation to optimality can be processed greedily 
 assigning the production levels to the cheapest plants for all scenarios.
 Their second  method uses a flow network relaxation which considers outage restrictions (CT 13)
as well as fuel consumption in an approximate fashion (CT7 to CT12) 
to deduce a tighter lower bound of the objective function than the simple greedy approach.
The computational effort to compute the bounds differs dramatically: while the greedy bounds are computed in less
than 10 seconds, solving the flow-network for all scenarios takes up to an hour for the biggest instances.
Their results are reported in this paper in Tables \ref{quickBounds} and \ref{bestDualROADEF}.

\section{MIP formulation without CT6 and CT12}\label{sec::v0MIP}

In this section, the constraints CT6 and CT12 are relaxed  similarly with \cite{lusby}.
It leads to a MIP formulation with  binary variables only for the outage decisions.

\begin{table}[ht]
      \centering
\caption{Set of variables}\label{setVar}
      \begin{tabular}{ll}
$d_{i,k,w}\in \{0,1\}$ & Beginning dates of outages decision.\\
$r_{i,k}\geqslant 0$ & Refueling levels. \\
$p_{i,k,s,t}\geqslant 0$ & Nuclear production of unit $i$ at cycle $k$ at $t$, 0 if $t$ is not in cycle $k$. \\
$p_{j,s,t}\geqslant 0$ & Nuclear production of unit $i$ at cycle $k$ at . \\
$r_{i,k}\geqslant 0$ & Refueling levels. \\
$x_{i,s}^{f}\geqslant 0$ & Fuel stock of T2 unit $i$ at the end of the optimizing horizon at scenario $s$. \\
$x_{i,k,s}^{init}\geqslant 0$ & Fuel levels at the beginning of production  cycle $k$ of unit $i$ . \\
$x_{i,k,s}^{fin}\geqslant 0$ & Fuel levels before the refueling $k+1$ of unit $i$, after production cycle $k$. \\
$x_{i,s,t}\geqslant 0$ & Fuel stock levels of T2 unit $i$ at time step $t$ and scenario $s$.
\end{tabular}

\end{table}

\subsection{Definition of variables}

\begin{figure}[ht]
      \centering
         \begin{minipage}[c]{.42\linewidth}
      \includegraphics[angle=0, width=.99\linewidth]{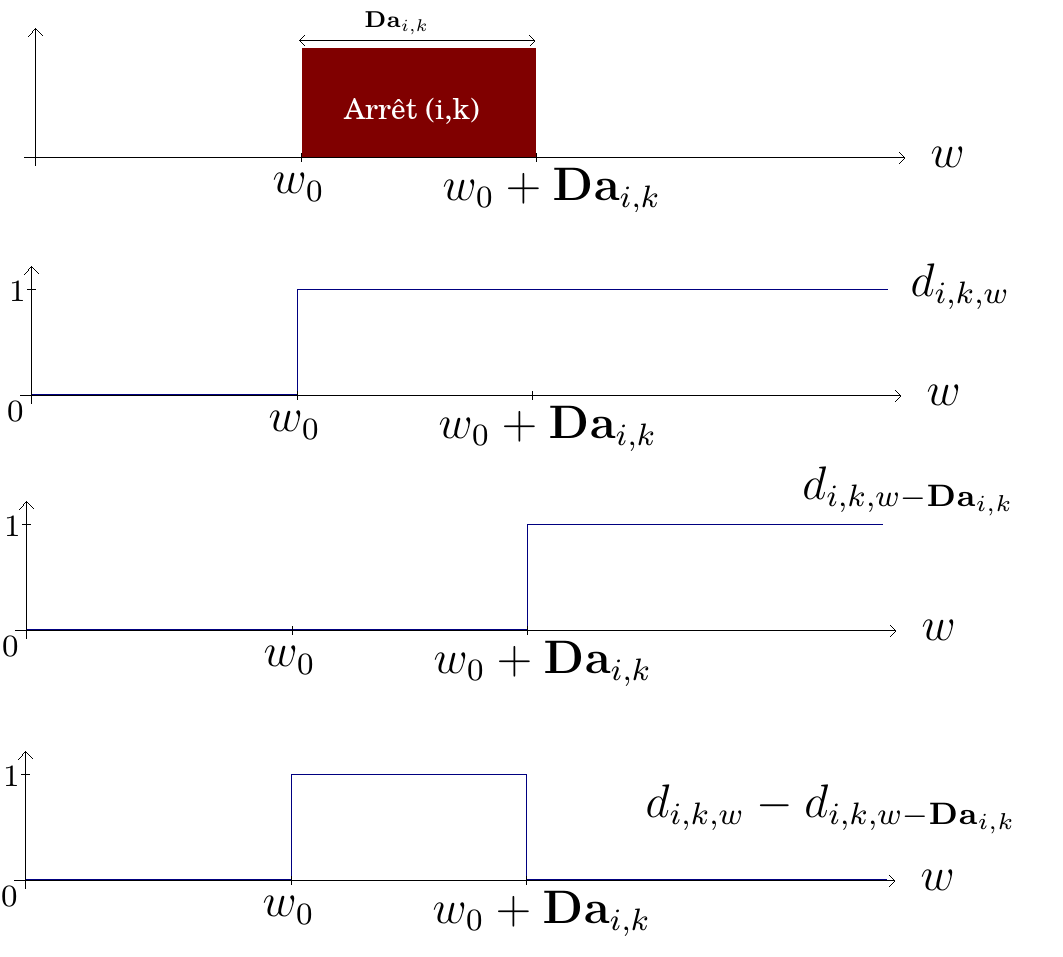}
	
   \end{minipage} \hfill
   \begin{minipage}[c]{.57\linewidth}
$\bullet$ $d_{i,k,w- {\mathbf{Da}_{i,k}}}$ indicates if  production cycle $(i,k)$  began before $w$.\\
$\bullet$ $d_{i,k,w}-d_{i,k,w-1} \in \{0,1\}$ is equal to  $1$ if and only if outage $(i,k)$ begins exactly at week  $w$.\\
$\bullet$  $d_{i,k,w- {\mathbf{Da}_{i,k}}}-d_{i,k,w-1- {\mathbf{Da}_{i,k}}} \in \{0,1\}$ is equal to  $1$ 
iff production campaign $(i,k)$ begins at week  $w$.\\
$\bullet$ $d_{i,k,w}-d_{i,k,w- {\mathbf{Da}_{i,k}}} \in \{0,1\}$ is equal to $1$ only during the outage of cycle $(i,k)$.\\
$\bullet$ $d_{i,k,w}-d_{i,k+1,w} \in \{0,1\}$ is equal to  $1$ only on the cycle $k$ of unit $i$.\\
$\bullet$ $d_{i,k,w- {\mathbf{Da}_{i,k}}}-d_{i,k+1,w} \in \{0,1\}$ is equal to  $1$ only on the production campaign $k$ of unit $i$.

   \end{minipage}

      \caption{Definition as "step variables" of binaries $d_{i,k,w}$ and related linear expressions}\label{DefVarD} 
\end{figure}

%
The binaries  $x_{i,k,w}$ in \cite{lusby} are equal to $1$ if and only if the beginning week for cycle $(i,k)$ is exactly $w$.
 Similarly with \cite{ucp}, we define the binaries  $d_{i,k,w}$ as ``step variables'' with 
 $d_{i,k,w}=1$  if and only if the outage beginning week for the cycle $k$ of the unit $i$ is before the week $w$, as illustrated Figure \ref{DefVarD}.
 We extend the notations with  $d_{i,k,w} =  0$ for $k > K_i$,$d_{i,-1,w} =  1$ for $w<0$,
  $d_{i,k,w} =  0$ for $w<0$ and $k>-1$.
  CT13 constraints, imposing that maintenance $(i,k)$ begins between $\mathbf{To}_{i,k}$ and $\mathbf{Ta}_{i,k}$
  reduce the definition of variables with  $d_{i,k,w}=0$ for $w < \mathbf{Ta}_{i,k}$ and $d_{i,k,w}=1$ for $w \geqslant \mathbf{Ta}_{i,k}$.


The other variables have a  continuous domain : 
 refueling quantities $r_{i,k}$ for each outage $(i,k)$,
T2 power productions $p_{i,k,t}$, fuel stocks at the beginning of campaign $(i,k)$ (resp at the end)  $x_{i,k}^{init},x_{i,k}^{fin}$,
T1 power productions $p_{j,t}$, and fuel stock $x_{i,s}^{f}$ at the end of the optimizing horizon.
We note that T2 power productions $p_{i,k,t,s}$ are duplicated for all cycle $k$ to have a linear model, 
$p_{i,k,t,s}=0$ if $t$ is not included in the production cycle $k$.
These variables are gathered in Table \ref{setVar}.

\subsection{MIP formulation}
Relaxing constraints CT6 and CT12, the previous definition of the variables allows to have a MIP formulation:

\begin{eqnarray}
v_0 =\min & \displaystyle\sum_{i,k} \mathbf{C}^{r}_{i,k} r_{i,k}  +
\sum_{j,s,t} \pi_s \mathbf{C}^{p}_{j,s,t} \mathbf{D}^t \: p_{j,s,t} - \sum_{i,s} \pi_s \mathbf{C}_{i,s}^{f} x_{i,s}^{f} 
 & \\
\forall i , k, w, & d_{i,k,w-1}\leqslant d_{i,k,w} & \label{PANprecedence}\\
\forall i , k, & d_{i,k,\mathbf{To}_{i,k}-1}\leqslant 0 & \label{PANtw0} \\
\forall i , k, & d_{i,k,\mathbf{Ta}_{i,k}}\geqslant 1 & \label{PANtw1}\\
\forall s,t, & \sum_{i,k} p_{i,k,s,t} + \sum_j p_{j,s,t} = \mathbf{Dem}^{t,s}\label{PANdemand}\\
\forall  j,s,t, & \mathbf{\underline{P}}_{j,t}^s \leqslant p_{j,s,t} \leqslant \mathbf{\overline{P}}_{j,t}^s\label{PANflexPower}\\
\forall i,k,s,t,  &  p_{i,k,s,t} \leqslant \mathbf{\overline{P}}_{i,t} (d_{i,k,w_t -{\mathbf{Da}_{i,k}}} - d_{i,k+1,w_t} )\label{PANcoupling}\\
\forall i,k, & \mathbf{\underline{R}}_{i,k} \; d_{i,k,W} \leqslant r_{i,k} \leqslant \mathbf{\overline{R}}_{i,k}\; d_{i,k,W}\label{PANrefuel}\\ 
\forall i,s, & x_{i,0,s}^{init} = \mathbf{Xi}_{i}\label{PANfuelInit}\\
\forall i,k,s, & x_{i,k,s}^{fin} = x_{i,k,s}^{init} - \sum_t \mathbf{D}^t \:  p_{i,k,s,t} \label{PANconso}\\
\forall i,k,s,  &  x_{i,k,s}^{init} - \mathbf{Bo}_{i,k} = r_{i,k} + \frac{\mathbf{Q}_{i,k} -1}{\mathbf{Q}_{i,k}} (x_{i,k-1}^{fin} - \mathbf{Bo}_{i,k-1})\label{PANpertes}\\
\forall i,k,s, &   x_{i,k,s}^{init}   \leqslant \mathbf{S}_{i,k} \label{PANmaxStock}\\
\forall i,k,s, & x_{i,k,s}^{fin} \leqslant \mathbf{A}_{i,k+1} + (\mathbf{S}_{i,k} - \mathbf{A}_{i,k+1})  (1-d_{i,k+1,W}) \label{PANanticip}\\
\forall i,k,s,  & x_{i,s}^{f} \leqslant x_{i,k,s}^{fin} + \overline{S}_{i}    ( 1-d_{i,k,W} + d_{i,k+1,W} )\label{PANfuelFinal}\\
  \forall c,w, & \sum_{(i,k)\in { \mathbf{A}}^c} (\alpha_{i,k,w} d_{i,k,w}) \leqslant \beta_w^c \label{PANordo}\\
  & d \in \{0,1\}^N, r,p,x \geqslant 0
\end{eqnarray}

Constraints (\ref{PANprecedence}) are required with definition of variables $d$.
Constraints (\ref{PANtw0}) and (\ref{PANtw1}) model CT13 time windows constraints: outage $(i,k)$ is operated between weeks $\mathbf{To}_{i,k}$ and $\mathbf{Ta}_{i,k}$.
Constraints (\ref{PANdemand}) model CT1 demand constraints.
Constraints (\ref{PANflexPower}) model CT2 bounds on T1 production.
Constraints (\ref{PANcoupling}) model CT3, CT4 and CT5 bounds on T2 production.
Constraints (\ref{PANrefuel}) model CT7 refueling bounds, with a null refueling when outage $i,k$ is not operated, ie $ d_{i,k,W}=0$.
Constraints (\ref{PANfuelInit}) write CT8 initial fuel stock.
Constraints (\ref{PANconso}) write CT9  fuel consumption constraints on stock variables of cycles $k$ $x_{i,k,s}^{init},x_{i,k,s}^{fin}$.
Constraints (\ref{PANpertes}) model CT10  fuel losses at refueling.
Constraints (\ref{PANmaxStock}) write CT11 bounds on fuel stock levels only on  variables $x_{i,k,s}^{init}$ which are the maximal stocks level over cycles $k$.
thanks to (\ref{PANconso}).
Constraints (\ref{PANanticip}) model CT11 minimum fuel stock before refueling, these constraints are active for a cycle $k$ only if the cycle is finished at the end of the optimizing
horizon, ie if $d_{i,k+1,W}=1$, which enforces to have disjunctive constraints where case $d_{i,k+1,W}=0$ implies a trivial constraints thanks to (\ref{PANmaxStock}).
Constraints (\ref{PANfuelFinal}) are linearization constraints to enforce $x_{i,s}^{f}$ to be the fuel stock at the end of the time horizon.
$x_{i,s}^{f}$ is indeed the $ x_{i,k,s}^{fin}$ such that $d_{i,k,W}=1$ and $d_{i,k+1,W}=0$, for the disjunctive constraints (\ref{PANfuelFinal}) 
that write a trivial constraints in the other cases thanks to (\ref{PANmaxStock}), we define $\overline{S}_{i}= \max_{k} \mathbf{S}_{i,k}$.
Constraints (\ref{PANordo}) are a common framework for scheduling constraints from CT14 to CT21, which was noticed independently in \cite{Jonc10,jost,lusby}. 
Distinguishing the cases, these constraints are written as following:


\begin{small}
\begin{eqnarray}
\forall c\in \CC_{14}, w, & \displaystyle
 \sum_{(i,k)\in { \mathbf{A14}^{c_{14}}}} (d_{i,k,w} -d_{i,k,w- ({ \mathbf{Da}_{i,k}} +  \mathbf{Se14}^{c_{14}})^+}) \leqslant 1 \label{CT14}\\
 \forall c\in \CC_{15}, w \in W_c^{15}, &  \displaystyle
 \sum_{(i,k)\in { \mathbf{A15}^{c}}} (d_{i,k,w} -d_{i,k,w- ({ \mathbf{Da}_{i,k}} +  \mathbf{Se15}^{c})^+}) \leqslant 1 \label{CT15}\\
 \forall c\in \CC_{16}, w, & \displaystyle \sum_{(i,k)\in { \mathbf{A16}^{c}}} (d_{i,k,w} -d_{i,k,w- { \mathbf{Se16}^{c}}}) \leqslant 1 \label{CT16}\\
 \forall c\in \CC_{17}, w \in \WW ,& \displaystyle \sum_{(i,k)\in { \mathbf{A17}^{c}}} (d_{i,k,w- { \mathbf{Da}_{i,k}}}
 -d_{i,k,w-{ \mathbf{Da}_{i,k}}- { \mathbf{Se17}^{c_{17}}}}) \leqslant 1 \label{CT17}\\
 \forall c\in \CC_{18}, w, & \displaystyle\sum_{(i,k)\in { \mathbf{A}_{18}^{c}}} (d_{i,k,w} -d_{i,k,w- { \mathbf{Se18}^{c}}} + d_{i,k,w- { \mathbf{Da}_{i,k}}}
 -d_{i,k,w- { \mathbf{Da}_{i,k}}- { \mathbf{Se18}^{c}}}) \leqslant 1 \label{CT18}\\
 \forall c_{19}, w, & \displaystyle\sum_{(i,k)\in { \mathbf{A19}}} (d_{i,k,w - {\mathbf{ L19}_{i,k}^{c_{19}}}}
 -d_{i,k,w - {\mathbf{ L19}_{i,k}^{c_{19}}} - { \mathbf{Tu19}_{i,k}^{c_{19}}}}) \leqslant \mathbf{Q19}^{c_{19}} \label{CT19}\\
 \forall c_{20}, w, &\displaystyle \sum_{(i,k)\in { \mathbf{A20}^{c_{20}}_w}} \: (d_{i,k,w} -d_{i,k,w- {\mathbf{Da}_{i,k}}}) \leqslant \mathbf{N20}_w^{c_{20}} \label{CT20}\\
 \forall c_{21}, w, & \displaystyle\sum_{i,k} \left (\sum_{t: w_t=w}\mathbf{\overline{P}}_i^t \right) (d_{i,k,w} -d_{i,k,w- {\mathbf{Da}_{i,k}}}) \leqslant \mathbf{Imax}_w^{c_{21}}\label{CT21}
\end{eqnarray}
\end{small}

$v_0$ is a dual bound for the whole problem, relaxing only constraints CT6 and CT12, requiring less binary variables than the formulation of \cite{Jonc10}.
Each dual bound proven for $v_0$ is thus a dual bound for the EURO/ROADEF 2010  Challenge.

\section{A lighter MIP formulation for CT6 constraints}

\begin{figure}[ht]
      \centering
      \includegraphics[angle=0, width=.67\linewidth]{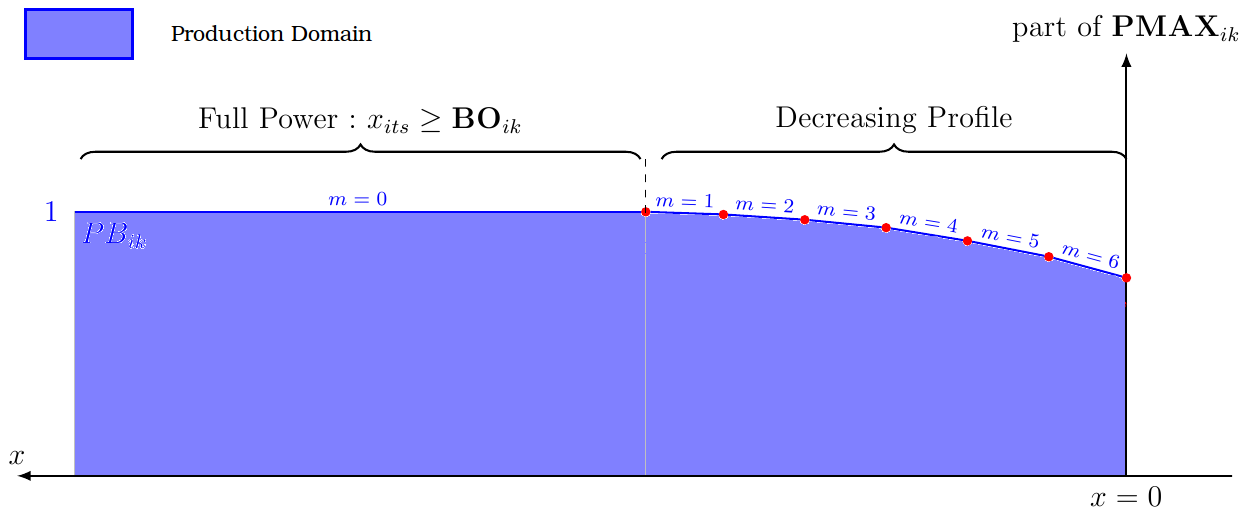}
      \caption{Illustration of production domain with light CT6 constraints}\label{evo1}  
\end{figure}

In the previous section,  CT6 constraints are relaxed, considering no  ``stretch'' decreasing profile as illustrated  Figure \ref{prodRoadef}.
In  \cite{Jonc10}, a exact compact formulation is provided, introducing binaries for all time steps, cycle and production mode $m \in \MM$,
which is not reasonable for the B\&B algorithm.
This section provides a lighter formulation for stretch constraints 
enforcing only upper bounds on the production as illustrated Figure \ref{evo1}.
This  is justified  because the T2 production has lower marginal costs than the T1 production.
The optimization tend to have T2 units producing at their maximal power,
 the upper bounds on the T2 production are essential.

As illustrated in Figure \ref{evo1}, the equation linking modes $m-1$ et $m$ is:
$$Y = \frac{\mathbf{c}_{i,k,m-1} - \mathbf{c}_{i,k,m}}{\mathbf{f}_{i,k,m-1} - \mathbf{f}_{i,k,m}} (X - \mathbf{f}_{i,k,m}) + \mathbf{c}_{i,k,m}$$
where $Y\in [0,1]$ is the ratio maximal power in stretch / maximal power $\mathbf{\overline{P}}_i^w$,
and $X$ denotes the residual fuel stock. 

To write mathematically  constraints to have upper bounds for the T2 production, we use that stretch decreasing profile  is concave.
Thus, the production domains are defined with the  intersection of the semi spaces defined with the equations of the different modes.
Adding in the model variables $x_{i,s,t}$ denoting the residual fuel stock of the unit $i$ at  time peroid $t$ for the scenario $s$,
the intersection of the semi spaces gives rise to following stretch constraints
\begin{equation}\label{ctStretch1}
\forall  i,k,s,t,m>0, \; \; \;  \frac {p_{i,k,s,t}} {\mathbf{\overline{P}}_i^t} \leqslant  
\frac{\mathbf{c}_{i,k,m-1} - \mathbf{c}_{i,k,m}}{\mathbf{f}_{i,k,m-1} - \mathbf{f}_{i,k,m}} (x_{i,s,t} - \mathbf{f}_{i,k,m}) + \mathbf{c}_{i,k,m}
\end{equation}



The definition of new variables $x_{i,s,t}$ requires following linking constraints to
enforce $x_{i,s,t}$ to be the residual fuel:
\begin{equation}\label{defVarStretch}
 \forall  i,k,s,t,\; \; \; x_{i,s,t} \leqslant x_{i,k,s}^{init} -  \sum_{t'\leqslant t} \mathbf{D}^{t'} \:  p_{i,k,s,t'}
+ M_i \;  (1 - d_{i,k,w_t} + d_{i,k-1,w_t}) 
\end{equation}
where $M_i = \max_k \mathbf{\overline{S}}_{i,k}$, such $M_i$  verifies $x_{i,t} \leqslant  M_i$. 
Indeed, if $d_{i,k,w} - d_{i,k-1,w}=1$, week $w$ happens in cycle $k$, the active constraint is
 $x_{i,t} \leqslant x_{i,k,s}^{init} -  \sum_{t'<t} \mathbf{D}^{t'} \:  p_{i,k,s,t'}$,
otherwise we have  $x_{i,t} \leqslant  M_i$, which is trivial thanks to the definition of $M_i$.
%

Adding constraints (\ref{ctStretch1}-\ref{defVarStretch}) in the MIP formulation of Section 4, it improved the dual bounds for the 
EURO/ROADEF 2010 Challenge. 

A special case is interesting: if  (  $\mathbf{c}_{i,k,m},\mathbf{f}_{i,k,m-1}$) do not depend on indexes  $k$,
we  can have a MIP formulation using less constraints than previously, 
writing the  constraints for the global production power $\sum_k {p_{i,k,w}}$:
\begin{equation}\label{ctStretch2}
\forall  i,s,t,m>0 \;\;\; \sum_k \frac {p_{i,k,s,t}} {\mathbf{\overline{P}}_i^t} \leqslant  \frac{\mathbf{c}_{i,m-1} - \mathbf{c}_{i,m}}{\mathbf{f}_{i,m-1} - \mathbf{f}_{i,m}} (x_{i,s,t} - \mathbf{f}_{i,m}) + \mathbf{c}_{i,m}
\end{equation}

\vskip 0.5cm
\noindent{\textbf{Remark}}: Adding light CT6 constraints with (\ref{ctStretch1}-\ref{defVarStretch}) or with (\ref{defVarStretch}-\ref{ctStretch2})
has few implications in terms of number of variables, adding $I \times S \times T$ continuous  variables   $x_{i,s,t}$,
whereas $I \times S \times T \times K$ continuous  variables were already in the model with T2 productions $p_{i,k,s,t}$.
The main difference appears in the number of constraints, there were mainly $I \times S \times T \times K$ constraints in the MIP of Section 4 with constraints (\ref{PANcoupling}),
 (\ref{defVarStretch}-\ref{ctStretch2}) require to add  $I \times S \times T \times (K+M)$ constraints,
 whereas (\ref{ctStretch1}-\ref{defVarStretch}) require to add  $I \times S \times T \times (M+1) \times K$ constraints.

\section{Parametric MIP relaxations}\label{sec::simpMIP}

This section propose new relaxations to compute dual bounds reducing significantly the size of the MIP problems.
A crucial point will be to prove that the restricted MIP computations induce dual bounds of $v_0$.

\subsection{Relaxations of all the outages }

A first simplification is the relaxation of outage constraints CT3.
Therefore, the null productions phases are relaxed, refueling are modeled only to compute the production costs of the T2 units.
The following MIP formulation expresses such problem where 
$d_{i,k}$ are the only binary variables corresponding to previous $d_{i,k,W}$ variables,
$x_{i,s}^{f}$ are the only stock variables required
and nuclear productions $p_{i,s,t}$ are not duplicated for all cycles $k$:

\begin{eqnarray}
v_3 =  \min & \displaystyle \sum_{i,k} \pi_s \mathbf{C}^{rld}_{i,k} r_{i,k} + 
\sum_{s,j,t} \pi_s  \mathbf{C}^{prd}_{j,t}  \mathbf{D}^t \: P_{j,s,t} - \sum_{i,s} \pi_s   \mathbf{C}_i^{val} x_{i,s}^{f}& \\
\forall i , k,    & d_{i,k+1} \leqslant d_{i,k}\\
\forall  i,k, & \mathbf{\underline{R}}_{i,k} \; d_{i,k} \leqslant r_{i,k}   \leqslant \mathbf{\overline{R}}_{i,k}\; d_{i,k}\\ 
\forall  t,j,s, & \mathbf{\underline{P}}_{j,t}^s \leqslant p_{j,s,t} \leqslant  \mathbf{\overline{P}}_{j,t}^s  \\
\forall t,s, & \sum_{i} p_{i,s,t} + \sum_j p_{j,s,t} = \mathbf{Dem}^{t,s}\\
\forall  i,s, & x_{i,s}^{f} \leqslant \Delta_i + \mathbf{Xi}_{i} + \sum_{k}  r_{i,k}   -  \sum_{t} \mathbf{D}^{t} \:  p_{i,s,t} \label{bornesAprouver}\\
\forall  i,k,s, & x_{i,s}^{f} \leqslant \mathbf{\overline{S}}_{i,k}\:  (d_{i,k} - d_{i,k+1}) \label{stockMaxFin}
\end{eqnarray}

Quantities $\Delta_i$ are defined for all T2 unit $i$ with $\Delta_i = \max_k (\mathbf{Bo}_{i,k}- \mathbf{Bo}_{i,0})$.
For all $(i,k)$,  $\mathbf{Bo}_{i,k}- \mathbf{Bo}_{i,0} \leqslant \Delta_i$
\vskip 0.3cm
\begin{proposition}  $v_3$ induced by the relaxation of all the outage gives dual bound to the Challenge ROADEF, with  
$v_3 \leqslant v_0$.
\end{proposition}

\textbf{Proof}: Let $(d_{i,k,w}^{*}, r_{i,k}^{*},p_{i,k,s,t}^{*}, p_{j,s,t}^{*},x^{*})$ an optimal solution of the  MIP defining $v_0$. 
Let us prove that $(d_{i,k,W}^{*}, r_{i,k}^{*},\sum_k p_{i,k,s,t}^{*},p_{j,s,t}^{*},x^{f*})$ is a feasible solution of the MIP defining $v_3$.
Constraints (\ref{bornesAprouver}) are the only constraints that are not trivially true.
Let $i \in \II, s\in\SS$.
Let $k$ the cycle , such that $x_{i,s}^{f*} = x_{i,k}^{fin*}$.
With $\frac{\mathbf{Q}_{i,k} -1}{\mathbf{Q}_{i,k}}\leqslant 1$,   (\ref{PANconso}) and (\ref{PANpertes}) imply : 
%
%
\small{
\begin{eqnarray*}
x_{i,k,s}^{fin*} + \sum_t \mathbf{D}^t  p_{i,k,s,t} = x_{i,k}^{init*} &  \leqslant & \mathbf{Bo}_{i,k}+  r_{i,k}^{*} + x_{i,k-1}^{fin*} - \mathbf{Bo}_{i,k-1}\\ 
x_{i,k-1,s}^{fin*} + \sum_t \mathbf{D}^t  p_{i,k-1,s,t} = x_{i,k-1}^{init*}  & \leqslant &   \mathbf{Bo}_{i,k-1}+  r_{i,k-1}^{*} + x_{i,k-2}^{fin*} - \mathbf{Bo}_{i,k-2}\\
\vdots & \vdots &  \vdots\\
x_{i,1,s}^{fin*} + \sum_{k,t} \mathbf{D}^t  p_{i,k,s,t} = x_{i,0}^{init*}  & \leqslant & \mathbf{Bo}_{i,0}+  r_{i,0}^{*} + x_{i,0}^{fin*}
 \end{eqnarray*}
}
\normalsize{
Adding these inequalities: $\displaystyle x_{i,s}^{f*} = x_{i,k,s}^{fin*}   \leqslant  \mathbf{Bo}_{i,k}- \mathbf{Bo}_{i,0} + r_{i,k} + x_{i,0,s}^{fin*} $.\\
It implies (\ref{bornesAprouver}), as we have
$x_{i,0,s}^{init*} = \mathbf{Xi}_{i}$ and $ \mathbf{Bo}_{i,k_0}- \mathbf{Bo}_{i,0} \leqslant \Delta_i$.\\
$(d_{i,k,W}^{*}, r_{i,k}^{*},\sum_k p_{i,k,s,t}^{*},p_{j,s,t}^{*},x^{fin \phantom{1} *})$ is thus a feasible solution of  the MIP defining $v_3$.
The cost of this solution is $v_0$, which is upper than the optimal solution, which proves $v_3 \leqslant v_0$.
$\square$
}

\subsection{Parametric relaxation of outages }

To improve the last relaxation, a parametric formulation   relaxes only the outages with an index $k>k^0$
and uses the MIP formulation of Section 4  for cycles $k\leqslant k^0$.
Hence, the positivity constraints of final stocks of production cycles applies for cycles $k< k^0$.
This gives rise to the following MIP formulation where constraints $M^{ordo}_{k^0} d \geqslant b_{k^0}^{ordo}$ 
gather the truncated constraints (\ref{CT14}-\ref{CT21})  considering only variables $d_{i,k,w}$ with index $k\leqslant k^0$.
As previously, we define quantities $\Delta'_i = \max_{k>k^0} (\mathbf{Bo}_{i,k_i^m}- \mathbf{Bo}_{i,k^0})$ to ensure that
for all $k>k^0$, $\mathbf{Bo}_{i,k_i^m}- \mathbf{Bo}_{i,k^0} \leqslant \Delta'_i$.

\begin{eqnarray}
v_3(k^0) =  \min & \displaystyle \sum_{i,k} \mathbf{C}^{rld}_{i,k} r_{i,k} + 
\sum_{j,s,t} \pi_s  \mathbf{C}^{prd}_{j,w}  \overline{D}^w \: p_{j,s,t} - \sum_{i,s}\pi_s   \mathbf{C}_i^{val} x_{i,s}^{f}& \\
\forall i , k,    & d_{i,k+1} \leqslant d_{i,k}\\
\forall i , k\leqslant k^0,    & d_{i,k} = d_{i,k,W}\\
& M^{ordo}_{k^0} d \geqslant b_{k^0}^{ordo} & \\
\forall  i,k, & \mathbf{\underline{R}}_{i,k} \; d_{i,k} \leqslant r_{i,k}   \leqslant \mathbf{\overline{R}}_{i,k}\; d_{i,k}\\ 
\forall i,s, & x_{i,0,s}^{init} = \mathbf{Xi}_{i} \\
\forall i,s,k\leqslant k^0 &   x_{i,k,s}^{init}   \leqslant \mathbf{\overline{S}}_{i,k} \\
\forall i,s,k\leqslant k^0,w& x_{i,k,s}^{fin} = x_{i,k,s}^{init} - \sum_t \mathbf{D}^t \:  p_{i,k,s,t} \\
\forall i,s,k< k^0, & x_{i,k,s}^{fin} \leqslant \mathbf{\overline{A}}_{i,k+1} + (\mathbf{\overline{S}}_{i,k} - \mathbf{\overline{A}}_{i,k+1}) \; (1-d_{i,k+1,W}) \\
\forall i,s,k\leqslant k^0  &  x_{i,k,s}^{init} - \mathbf{Bo}_{i,k} = r_{i,k} + \frac{\mathbf{Q}_{i,k} -1}{\mathbf{Q}_{i,k}} (x_{i,k-1,s}^{fin} - \mathbf{Bo}_{i,k-1}) \\
\forall i,s,k\leqslant k^0  & x_{i,s}^{fin} \leqslant x_{i,k,s}^{fin} + \mathbf{\overline{S}}_{i}    (d_{i,k,W} - d_{i,k+1,W} )\\
\forall  t,j,s, & \mathbf{\underline{P}}_{j,t}^s \leqslant p_{j,s,t} \leqslant  \mathbf{\overline{P}}_{j,t}^s  \\
\forall t,i,k\leqslant k^0,s,  &  p_{i,k,s,t} \leqslant \mathbf{\overline{P}}_{j,t}^s (d_{i,k,w_t -{\mathbf{Da}_{i,k}}} - d_{i,k+1,w_t} )\label{bornesAprouver02}\\
\forall t,s, & \sum_{i,k\leqslant k^0} p_{i,k,s,t} + \sum_j p_{j,s,t} = \mathbf{Dem}^{t,s}\label{agregPowCycle}\\
\forall  i,k\leqslant k^0,s, & x_{i,s}^{f} \leqslant \ x_{i,k,s}^{init}   -  \sum_{t} \mathbf{D}^{t} \:  p_{i,k,t} + \mathbf{\overline{S}}_{i,k}\:  (1 + d_{i,k+1} - d_{i,k})\\
\forall  i,s, & x_{i,s}^{f} \leqslant \Delta'_{i} + x_{i,k^0,s}^{init}  + \sum_{k>k^0}  r_{i,k}   -  \sum_{t} \mathbf{D}^{t} \:  p_{i, k^0, t} + \mathbf{\overline{S}}_{i,k}\:  (1 - d_{i,k^0} )\label{bornesAprouver2}\\
 & d \in \{0,1\}^N, r,p,x \geqslant 0
\end{eqnarray}

\begin{proposition}  
For all  $k^0\in \KK$, 
$v_3(k^0) \leqslant v_0$.
Each dual bound proven for a $v_3(k^0)$ is thus a dual bound for the EURO/ROADEF 2010  Challenge.

\end{proposition}

\noindent{\textbf{Proof}}: Let $(d_{i,k,w}^{*}, r_{i,k}^{*},p_{i,k,s,t}^{*}, p_{j,s,t}^{*},x^{*})$ be an optimal solution of the MIP defining $v_0$. 
Like previously,  we prove first that we have a feasible solution of the MIP defining $v_3(k^0)$
with $d_{i,k} =d_{i,k,W}^{*}$, $d_{i,k,w} =d_{i,k,w}^{*}$ for $k\leqslant k^0$, $p_{j,s,t} =p_{j,s,t}^{*}$ ,  $p_{i,k,s,t} =p_{i,k,s,t}^{*}$ for $k< k^0$
$p_{i,k^0,s,t} = \sum_{k\geqslant k^0} p_{i,k,s,t}^{*}$ and $x=x^{*}$. 
 (\ref{bornesAprouver02}) and (\ref{bornesAprouver2}) are the only constraints that are not trivially verified.

$0 \leqslant p_{i,k,s,t}^{*} \leqslant \overline{\overline{P}}_{j,t}^s(d_{i,k,w_t -{\mathbf{Da}_{i,k}}} - d_{i,k+1,w_t})$. 
It implies:\\ 
$p_{i,k^0,s,t} = \sum_{k\geqslant k^0} p_{i,k,s,t}^{*} \leqslant \overline{\overline{P}}_{j,t}^s \sum_{k\geqslant k^0}(d_{i,k,w_t -{\mathbf{Da}_{i,k}}} - d_{i,k+1,w_t-\mathbf{Da}_{i,k+1}})$,\\
using $d_{i,k+1,w -\mathbf{Da}_{i,k+1}} \leqslant d_{i,k+1,w}$. 
The telescopic summation and $d_{i,K+1,w_t} = 0$ imply (\ref{bornesAprouver02}). 
%
%
To prove (\ref{bornesAprouver2}), we use $k_i^m$ the last cycle operated in the solution $(d_{i,k,w}^{*})$, i.e. 
$k_i^m$ is the maximal $k$ such that $d_{i,k,W}^{*}\neq 0$. 
(\ref{bornesAprouver2}) are trivially verified if $d_{i,k^0}=1$ with (\ref{stockMaxFin}), i.e. if $k_i^m< k^0$.
We suppose $k_i^m \geqslant k^0$.
With  $\frac{\mathbf{Q}_{i,k} -1}{\mathbf{Q}_{i,k}}\leqslant 1$ and (\ref{PANpertes}), 
we have for all $k>0$, 
$x_{i,k,s}^{fin*} + \sum_t \overline{D}^w  p_{i,k,s,t} = x_{i,k,s}^{init*}  \leqslant \mathbf{Bo}_{i,k}+  r_{i,k}^{*} + (x_{i,k,s}^{fin*} - \mathbf{Bo}_{i,k-1})$.
Adding these inequalities from $k^0+1$ to $k_i^m$:\\
$\displaystyle x_{i,k_i^m,s}^{fin}  + \sum_{k=k^0}^{k_i^m}\sum_{t} \overline{D}^w  p_{i,k,s,t}  \leqslant  \mathbf{Bo}_{i,k_i^m}- \mathbf{Bo}_{i,k^0} + r_{i,k} + (x_{i,k^0,s}^{fin} - \mathbf{Bo}_{i,k^0})$\\
These last inequalities imply (\ref{bornesAprouver2}), as
we have $x_{i,k^m,s}^{fin}=  x_{i,s}^{f}$,  $p_{i,k,s,t}=0$ and $r_{i,k}=0$ for all $k>k_i^m$ and $ \mathbf{Bo}_{i,k_i^m}- \mathbf{Bo}_{i,k^0} \leqslant \Delta_i$.


$(d_{i,k,W}^{*}, r_{i,k}^{*},\sum_k p_{i,k,s,t}^{*},p_{j,s,t}^{*},x^{fin \phantom{1} *})$ is thus a feasible solution of  the MIP defining $v_3(k^0)$.
The cost of this solution is $v_0$, which is upper than the optimal solution, which proves 
$v_3(k^0) \leqslant v_0$.
$\square$

\section{Dual bounds with aggregations and problem reductions}

In this section, we prove that dual bounds for the challenge ROADEF can be calculated with three reductions:
 exact preprocessing,  aggregation of  time steps and also
reduction  of the number of  scenarios.
These reductions apply to reduce any MIP giving dual bounds for the Challenge ROADEF,
resulting from sections 4, 5 or 6.

%

\subsection*{\textbf{7.1 Exact preprocessing}}

This section aims to reduce the size of the MIP computation to solve applying exact preprocessing to delete variables.
 Two propositions are mentioned here do deal with less variables in the MIP computations, the proofs (not difficult) are given in the Appendix A.
First,  Proposition \ref{tightenTW} allows to reduce time windows by propagation.
Then, Proposition \ref{tightenContVar} allows to deal with less continuous variables.

\begin{proposition}  \label{tightenTW}
Denoting
$\mathbf{Lmin}_{i,k} = \left\lceil \frac {\mathbf{Rmin}_{i,k} - \mathbf{Amax}_{i,k}} {\mathbf{D}^w  P_i}\right\rceil$ with 
$P_i = \max_w \mathbf{Pmax}_{i,w}$, 
$\mathbf{To}_{i,k}$ and $\mathbf{Ta}_{i,k}$ can be strengthened in $\widetilde{To}_{i,k}$ and $\widetilde{Ta}_{i,k}$ with induction relations:
\begin{equation}\label{tightenTWo}
\forall i \in \II, k>0, \;\;\; \widetilde{To}_{i,k} = \max(\mathbf{To}_{i,k}, \widetilde{To}_{i,k-1} + \mathbf{Da}_{i,k-1} + \mathbf{Lmin}_{i,k})
\end{equation}
\begin{equation}\label{tightenTWa}
\forall i \in \II, k<K, \;\;\; \widetilde{Ta}_{i,k} = \min(\mathbf{Ta}_{i,k}, \widetilde{To}_{i,k+1} - \mathbf{Da}_{i,k} - \mathbf{Lmin}_{i,k})
\end{equation}

$\widetilde{To}_{i,k}$ are first computed by induction with $k$ increasing.
Then $\widetilde{Ta}_{i,k}$ are computed by induction with  $k$ decreasing.
\end{proposition}

%

%
%


\begin{proposition} \label{tightenContVar}

Denoting for all $(m,n) \in K$, $\mathbf{q}_{m,n} = \prod_{l=m}^{n} \dfrac{\mathbf{Q}_{i,l} -1}{\mathbf{Q}_{i,l}}$,
we have following relations which allow to delete variables $ x_{i,k,s}^{init},  x_{i,k,s}^{fin}$:

\begin{scriptsize}
\begin{equation}
 x_{i,k,s}^{init} = \mathbf{q}_{1,k} \mathbf{Xi}_{i}
+ \sum_{l=0}^{k-1} \mathbf{q}_{l+2,k} 
\left(  r_{i,l+1} - \mathbf{q}_{l+1,l+1} \left( \sum_t \mathbf{D}^t \:  p_{i,l,s,t} - \mathbf{Bo}_{i,l} \right)   \right)
\end{equation}
\begin{equation}
x_{i,k,s}^{fin} = \mathbf{q}_{1,k} \mathbf{Xi}_{i}
+ \sum_{l=0}^{k-1} \mathbf{q}_{l+2,k} 
\left(  r_{i,l+1} - \mathbf{q}_{l+1,l+1} \left( \sum_t \mathbf{D}^t \:  p_{i,l,s,t} - \mathbf{Bo}_{i,l} \right)   \right) - \sum_t \mathbf{D}^t \:  p_{i,k,s,t} 
\end{equation} 
\end{scriptsize}


\end{proposition}

%



\subsection*{\textbf{7.2 Dual bounds by time-step aggregation}}


In this section, we prove that time step aggregation for the nuclear power plants allows to compute dual bounds.
Therefore, we denote the following MIP as a general expression for a MIP giving lower bounds for the Challenge ROADEF after section 4, 5 or 6, with or without exact preprocessing.


\begin{eqnarray}
v =  \min_{x\in X,y,z\geqslant 0} & \displaystyle {c}_x x + \sum_{t} \mathbf{D}^t {c}_{t} y_{t} & \\
s.t: 
 & T^1 x + \sum_{t} \mathbf{D}^t W^{1}  y_{t}  \geqslant h^1 \label{ctSlavAgg}\\
\forall t, & T_t^2 x +  W^{2}  y_{t}  \geqslant h^2_{t} \label{ctSlavTemp}
\end{eqnarray}

where $x$ denotes the vector of variables $d_{i,k,w}$, $y_t\geqslant 0$ consider only the T1 production variables
whereas $z_t\geqslant 0$ denotes the other continuous variables which do  not impact the cost.


An important point is that $ W^{1},W^{2}$ do not depend on $t$.
Two hypotheses, verified in the datasets B and X of the challenge, are essential in the following proofs.
First, $\mathbf{D}^t$ is constant over time, we denote $\mathbf{D} =\mathbf{D}^t$ in that context.
Then, T1 production costs are constant over weeks, defining $ c_w$ quantities:
\begin{equation}\label{hypCostConst}
\forall w \in \WW, \forall (t,t'), w_t = w_{t'} \Longrightarrow    c_{t} = c_{t'}= c_w
\end{equation}
 
\vskip 0.3cm
 
Let $\overline{D}^w$ the duration of a week in the time unit of $t$, $\displaystyle \overline{D}^w = \sum_{t, w_t=w} \mathbf{D}^t$,
which is common for all weeks.
Let
\begin{small}
$\displaystyle \overline{T}_{w}^2  = \sum_{t, w_t=w} \alpha  T_{t}^2$,
$\displaystyle \overline{h}^2_{w}  = \sum_{t, w_t=w} \alpha  h^2_{t}$.
\end{small}

\begin{eqnarray}
\overline{v}                                =  \min_{x\in X,y\geqslant 0} & \displaystyle {c}_x x + \sum_{w}\overline{D}^w \overline{c}_{w} y_{w} & \\
s.t:  
 & T^1 x + \sum_{w} \overline{D}^w {W}^1  y_{w} \geqslant h^1 \label{ctSlavAgg2}\\
\forall w, & \overline{T}_w^2  x +  {W}^2  y_{w} \geqslant \overline{h}^2_{w} \label{ctSlavTemp2}
\end{eqnarray}

\begin{proposition}\label{dualboundAgregTimeStep}
If the hypothesis \ref{hypCostConst} is valid,
aggregating production time steps to weeks provide a dual bound for the disaggregated problem: 
$ \overline{v} \leqslant {v} $.
\end{proposition}

\noindent{\textbf{Proof}}: Let $(x^{*},y^{*})$ 
an  optimal solution of the MIP  defining $v$. 
Let 
$\displaystyle \overline{y}_w^{*}  = \alpha \sum_{t, w_t=w}  y_t^{*}$, 
Let us prove that $(x^{*}, \overline{y}_w^{*})$ is feasible in the MIP defining $\overline{v}$. 
To prove  (\ref{ctSlavAgg2}), we notice that for all $w$:\\
$\displaystyle\sum_{t, w_t=w} \mathbf{D} W^1  y^*_{t} =  W^1 \sum_{t, w_t=w} \mathbf{D}   y^*_{t}
= W^1 \sum_{w} \mathbf{D}^w   \overline{y}^*_{w}$.\\
$\overline{y}_w^{*}$ verifies (\ref{ctSlavAgg2}) thanks to (\ref{ctSlavAgg}) and this last equality.\\
(\ref{ctSlavTemp2}) are proven from (\ref{ctSlavTemp})  aggregating constraints relatives to time steps $t$ for all week $w$ with   $w_t = w$ with weight $\alpha$:
$$ \sum_{t, w_t=w}\alpha T^2_{t} x^* +   \sum_{t, w_t=w} \alpha W^2  y_{t}^* \geqslant \sum_{t, w_t=w} \alpha h^2_{t} = \overline{h}^2_{w} $$

\noindent{With $\displaystyle \sum_{t, w_t=w} \alpha T^2_{t} x^* +   \sum_{t, w_t=w} \alpha W^2  y_{t}^* =
\alpha \left(\sum_{t, w_t=w} T^2_{t} \right)  x^* +   W^2  \left(\sum_{t, w_t=w} \alpha y_{t}^* \right)$,}

\noindent{$\displaystyle \sum_{t, w_t=w} \alpha T^2_{t} x^* +   \sum_{t, w_t=w} \alpha W^2  y_{t}^* = \overline{T}_w^2  x^* +  {W}^2  y_{w}^*$.}\\
The aggregated inequality above is thus exactly (\ref{ctSlavTemp2}).

%
The objective cost associated with $(x^{*}, \overline{y}_w)$ is thus superior to  the optimum $\overline{v}$:
$$\displaystyle \overline{v} \leqslant  {c}_x x^{*} + \sum_{w} \overline{D}^w \overline{c}_{w} \overline{y}_{w} = {c}_x x^{*} + \sum_{t} \overline{D}^t \overline{c}_{t} {y}^{*}_{t}$$ 
We note that last equality is true thanks to the hypothesis (\ref{hypCostConst}) $\square$.


\vskip 0.5cm
\noindent{\textbf{Remark}}: 
The hypothesis (\ref{hypCostConst}) is valid for the  datasets B and X of the Challenge, not for the dataset A.
The aggregation of production time steps to weekly time steps
  provides dual bounds of $v_0$ for datasets B and X, justifying that the Benders decomposition in \cite{lusby} allows to compute 
dual bounds for the Challenge ROADEF 2010. 
Their LP relaxation with Benders decomposition provide thus dual bounds, 
but also the dual bounds furnished after the Cut\&Branch phase.
However, our proof use hypotheses that does not guarantee that the time step aggregation furnish dual bounds for the dataset A.

\subsection*{\textbf{7.3 Dual bounds by scenario decomposition}}

In this section, we face another  bottleneck for an efficient solving: the number of scenarios.
Therefore, we define the following MIP as a general expression for the MIP giving bound after section 4, 5 or 6, with exact preprocessing
and possibly the lastly aggregation.

\begin{eqnarray}
v^{sto} =  \min_{x\in X,y\geqslant 0} & \displaystyle {c}_x x + \sum_{s}  \overline{c}_{s} y_{s} & \\
s.t:  & A x  \geqslant a \label{ctMaster03}\\
\forall s, & T x + {W}  y_{s} \geqslant h_s \label{ctSlavAgg03}
\end{eqnarray}

In this denomination, $x$ denotes the $x$ denotes the first stage variables $d_{i,k,w}, r_{i,k}$, $y_s\geqslant 0$ gather the other continuous variables
duplicated for all scenario $s$. Constraint matrices $A$, $W$ and $T$ do not depend from $s$.

We denote $v_{s}^{det}$ the following deterministic MIP for all scenarios $s \in\mathcal{S}$:

\begin{eqnarray}
v_{s}^{det} =  \min_{x\in X,y\geqslant 0} & \displaystyle {c}_x x +  {c}_{s} y & \\
s.t:  & A x  \geqslant a \label{ctMaster3}\\
 & T x +  {W}  y \geqslant h^1_s \label{ctSlavAgg3}
\end{eqnarray}

\begin{proposition} We have $\sum_s \pi_s v_{s}^{det} \leqslant v^{sto}$.
In other words, dual bounds for the whole ROADEF problem can be calculated with $|\SS|$ independent parallel computations of dual bounds
on  reduced problem with single scenarios. 
 \end{proposition}\label{prop::desagSto}

\noindent{\textbf{Proof}}: We reformulate the problem  duplicating first stage variables for all scenarios $x_s =x$ 
and using relation $\sum_{s\in\mathcal{S}} \pi_s =1$:\\

\begin{eqnarray*}
v^{sto}  =  \min_{x_s\in X,y\geqslant 0} & \displaystyle \sum_s \pi_s {c}_x x_s + \sum_{s} \pi_s {c}_{s} y_{s} & \\
s.t:  & A x  \geqslant a \\
 \forall s, & x^s=x \\
 \forall s, & A x_s  \geqslant a \\
 \forall s, & T x + {W}  y_{s} \geqslant h^1_s 
\end{eqnarray*}

Relaxing constraints $x^s =x$, we get a dual bound for $v_{1}$.
This relaxation implies independent sub-problems for all scenarios:

 \begin{eqnarray*}
v^{sto}  \geqslant \min & \displaystyle \sum_s \pi_s {c}_x x^s + \sum_{s} \pi_s {c}_{s} y_{s} & =  \displaystyle \sum_s \pi_s \left(\min_{x_s,y_s} {c}_x x^s + \sum_{s} {c}_{s} y_{s}\right)\\
s.t:  \forall s, & A x_s  \geqslant a & \hskip 1.5cm s.t:  \forall s,  A x_s  \geqslant a \\
 \forall s, & T x_s + \sum_{w} {W}  y_{s} \geqslant h^1_s & \hskip 1.5cm \forall s,  T x_s + \sum_{w} {W}  y_{s} \geqslant h^1_s 
 \end{eqnarray*}

 \noindent{The sub-problem of the scenario $s$ has value $ \pi_s v_{s}^{det}$, it proves $\sum_s \pi_s v_{s}^{det} \leqslant v^{sto}$. $\square$}

\vskip 0.3cm


One can notice that the previous proof can be generalized with sub-problems containing several scenarios. 
Let  $S^0 \subset \mathcal{S}$. Let ${v}_{S^0}$ the following restriction: 
\begin{eqnarray}
{v}_{S^0} =  \min_{x\in X,y\geqslant 0} & \sum_{s \in S^0} \pi_s \displaystyle {c}_x x + \sum_{ s \in S^0}  {c}_{s} y_{s} & \\
s.t:  & A x  \geqslant a \label{ctMaster4}\\
\forall s\in S^0, & T^1 x +  \overline{W}  y \geqslant h^1_s \label{ctSlavAgg4}
\end{eqnarray}


One can calculate dual bounds for the challenge ROADEF
with restricted computations containing a restricted number of  scenarios: 

\begin{proposition}  Let $S_n$ for $n \in \NN$ a partition of  $\mathcal{S}$, i.e. $\mathcal{S}= \bigcup_{n \in \NN} S_n$ with $S_n$ 
being disjoint subsets.
We have :
\begin{equation}
 \displaystyle \sum_s \pi_s v_{s}^{det} \leqslant \sum_{n \in \NN} {v}_{S_n}  \leqslant v^{sto} 
\end{equation}
 \end{proposition}

 \noindent{\textbf{Proof}}: The proof of ${v}_{S_n}$ is identical to the previous proof with relaxations $x^s =x^{s'}$ only when $s,s'$
 are not in the same subset $S_n$ of the chosen partition.
 Having this bounds, one can operate the relaxations $x^s=x$ in each sub-problem defined by a partition element like previously,
 which ensures $ \sum_{n \in \NN} {v}_{S_n}  \leqslant v^{sto} $.
 Then, for all subproblems ${v}_{S_n}$, we relax the remaining constraints $x^s=x$ like in the previous proposition, 
 which ensures $ \sum_s \pi_s v_{s}^{det} \leqslant \sum_{n \in \NN} {v}_{S_n} $. $\square$

 .

\section{Computational results}

This section analyzes the computational results.
We report firstly the characteristics of the datasets of the Challenge.
Before giving the dual bounds, the solving capabilities  of the B\&B algorithm are analyzed for the MIP formulations.
Our implementation used OPL and Cplex version 12.5 to solve MIP and LP problems.
Our experimentations were computed with a laptop running Linux Ubuntu 12.04 with an  Intel Core2 Duo processor, 2.80GHz.

\begin{table}
\centering
\caption{Characteristics of the instances of the 2010 EURO/ROADEF Challenge:
number of T1 and T2 power plants, number of production and weekly time steps and number of binaries.
nbVar0 is the total number of variables without preprocessing, total amplitude of the time windows.
nbVar2 is the total number of variables with the preprocessing of Proposition  \ref{tightenTW}. 
nbVar1 and nbVar3 are respectively the remaining variable after Cplex preprocessing
with respectively without and  with  preprocessing of Proposition  \ref{tightenTW}.
}\label{instancesROADEF}\label{instancesTronquees} 

\begin{tabular}{|l|c c c c|c c|c c |c c |}
\hline
Data&I&J &K&S&T&W&nbVar0 & nbVar1&nbVar2 & nbVar3\\

\hline
A1&10&11&6&10&1750&250&3892&463&483&424\\
A2&18&21&6&20&1750&250&7889& 961&892&761\\
A3&18&21&6&20&1750&250&8162&875&841&698\\
A4&30&31&6&30&1750&250&17465&875&841&698\\
A5&28&31&6&30&1750&250&15357&2797&2750&2494\\
\hline
B6&50&25&6&50&5817&277&24563&3466&3467&3054\\
B7&48&27&6&50&5565&265&35768&6435&9052&5846\\
B8&56&19&6&121&5817&277&69653&22482&30626&20763\\
B9&56&19&6&121&5817&277&69306&22482&30626&20763\\
B10&56&19&6&121&5565&265&29948&4236&5084&3790\\
\hline
X11&50&25&6&50&5817&277&20081&3478&3499&3216\\
X12&48&27&6&50&5523&263&27111&4348&5321&4035\\
X13&56&19&6&121&5817&277&30154&4697&4403&4104\\
X14&56&19&6&121&5817&277&30691&5378&6088&4879\\
X15&56&19&6&121&5523&263&27233&3992&4372&3618\\
\hline
\end{tabular}
\end{table}
%
%


\subsection{\textbf{Data characteristics}}

Three datasets were provided by EDF for the Challenge.
These three datasets are non-confidential and now available online. 
The data characteristics are provided in Table \ref{instancesROADEF}.
Dataset A contains five rather small instances given for the qualification phase of the Challenge.
Production time steps are discretized daily for the instances of dataset A, in a horizon of five years, with 10 to 30 nuclear units having $6$ 
production cycles, and 10 to 30 stochastic scenarios of demands.
Instances B and X  are more representative of real-world size instances.
Production time steps are discretized with $8h$ time steps to analyze the  impact of daily variability of power demands, in a horizon of five years, with 20 to 30 T1 units, around 50 nuclear units,
and 50 to 120 stochastic scenarios of demands.
These datasets are now public, dataset X was secret for the challenge, with instance characteristics that had to be similar to dataset B.
This symmetry holds for the number of units, cardinal of sets of Table \ref{setROADEF} and number of
constraints of each type,
but this does not hold for the number of binaries (i.e. the cumulated amplitude of time windows) required to solve  the instances: 
Table \ref{instancesROADEF} shows in the column nbVarBin the number of binary variables,.
Instances B8 and B9 are more combinatorial, with no time windows constraints for cycles $k\geqslant 3$.

%
%
%



\begin{table}
\centering
   \caption{MIP Convergence characteristics on $v_s^{det}$ applied to the mean scenario:
    comparison of  the gap of dual and primal bounds to the best primal solution known after LP relaxation, after Cplex cuts at the root node, 
    and in one hour MIP solving.
    Value are best primal solution known for $v_0$ computed with \cite{dupin2015modelisation}}\label{BandBSansStretch125}
\begin{tabular}{|l|l |l l |l l l|l l l|}
\hline
&Value&LP&CPU&inf0&sup0&CPU&inf&sup&CPU\\
\hline
A1&152170M&0,17\%&0,09&0,00\%&0,00\%&1,1&0,00\%&0,00\%&1,1\\
A2&145201M&0,22\%&0,4&0,02\%&0,00\%&3,8&0,00\%&0,00\%&4,9\\
A3&152582M&0,36\%&0,11&0,05\%&0,00\%&4,4&0,00\%&0,00\%&3\\
A4&102421M&0,95\%&3,1&0,42\%&0,09\%&21,4&0,00\%&0,00\%&334\\
A5&119785M&0,90\%&7,5&0,53\%&0,44\%&55&0,15\%&0,00\%&3600\\
\hline
B6&76966M&2,23\%&11,2&0,61\%&0,51\%&64&0,13\%&0,00\%&3600\\
B7&74233M&2,90\%&66,7&0,78\%&1,99\%&346&0,52\%&0,98\%&3600\\
B8&73239M&9,47\%&244&7,95\%&NS&3320&7,95\%&NS&3600\\
B9&72812M&8,51\%&350&6,69\%&NS&3600&6,69\%&NS&3600\\
B10&69501M&3,45\%&14,7&0,66\%&0,27\%&77&0,07\%&0,00\%&3600\\
\hline
X11&73018M&1,65\%&12,1&0,58\%&0,46\%&124&0,37\%&0,01\%&3600\\
X12&70604M&2,98\%&21,6&0,58\%&0,35\%&145&0,20\%&0,00\%&3600\\
X13&69230M&2,61\%&17&1,30\%&1,56\%&411&1,09\%&0,22\%&3600\\
X14&68395M&3,09\%&25&1,22\%&0,81\%&267&0,89\%&0,04\%&3600\\
X15&66028M&3,60\%&13,4&0,47\%&0,45\%&100&0,10\%&0,00\%&3600\\
\hline
\end{tabular}

\end{table}

\subsection{\textbf{MIP solving characteristics for single scenarios instances}}

\paragraph{\textbf{MIP relaxation of CT6 and CT12}}
First results justify the commonly used simplification to aggregate production time steps to weekly time steps, and to aggregate scenarios.
Using Cplex 12.3, the B\&B algorithm is completely inefficient on instances B8 et B9 even with time steps and scenario aggregations and relaxation of
difficult constraints CT6 and CT12, 1h resolution time is not enough to compute the LP relaxation.
The size of these instances is  a strong limiting factor.
For the other instances, the  B\&B  search is efficient, with low gaps between the best primal and dual bounds  in 1h, we refer to Table \ref{BandBSansStretch125}.
With Cplex 12.5, dual bounds can be computed for all instance with a single scenarios, as shown in Table \ref{BandBSansStretch125}.

 Furthermore, we observed that the variable definition is crucial in the efficiency of the branching quality:
our level variables $d_{i,k,w}\leqslant d_{i,k,w+1}$  imply  better branching and MIP convergence
than using binaries $x_{i,k,w}=d_{i,k,w}- d_{i,k,w-1}$ as in \cite{lusby} with GUB constraints $\sum_w x_{i,k,w}\leqslant 1$.
These results are coherent with \cite{ucp,vielma}.

We note that the exact preprocessing proposed in section 7.1 is shown useful in Table \ref{instancesROADEF}: 
our  tailored preprocessing of time windows of section 7.1 is not redundant with the one from Cplex,
which is particularly useful to deal with less binaries for the difficult instances B7, B8 and B9.

\paragraph{\textbf{MIP with light CT6 constraints}}
Adding upper bounds for CT6 stretch constraints, the computations of the LP relaxation are not possible for the datasets B and X,
even with the more compact formulation (\ref{ctStretch2}).
On dataset A, it increases on average  the number of variables and constraints before Cplex preprocessing with factors $1,5$ and $4,5$ respectively.
After Cplex preprocessing, it remains factors $1,3$ and $3,7$ which increases dramatically the LP relaxation consumption of CPU and memory. 
This point was a bad surprise : this formulation compared to \cite{Jonc10}
tried to get rid of too many binaries using only extra continuous variables, 
hoping it will not degrade the B\&B convergence characteristics.
The branching structure is indeed similar for small instances, but the difficulty is here in the size of the problems due to many continuous variables.
Furthermore, the gain in LP relaxation quality was measured on dataset A.
Another negative result was that these gains were not very significant (in the order of $0,1\%$).
In the compromise search between computation time and dual bound quality, this justifies
to relax constraints CT6 in the following.
Such results justify also the simplification hypotheses from  \cite{lusby}.

\begin{figure}[ht]
      \includegraphics[angle=0, width=.997\linewidth]{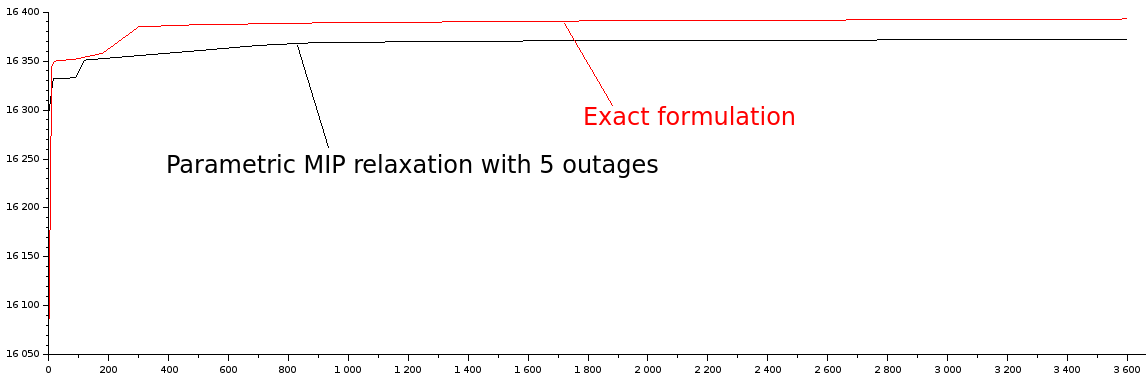}
      {for instance B6 }  
      \vskip 0.3cm
             \includegraphics[angle=0, width=.997\linewidth]{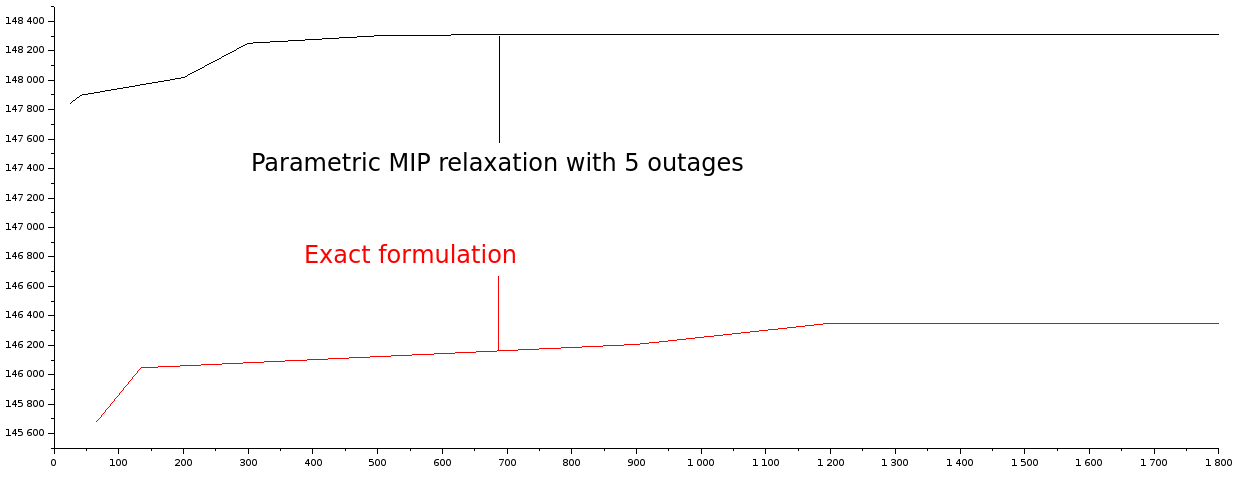}
             {for instance B7 }  
      \vskip 0.3cm
      \caption{Comparison of dual bound convergence for a deterministic MIP computation
      for the exact formulation without CT6 and CT12 and for the parametric relaxation with $k_0=5$  
      for instance B6 (relatively easy) and for difficult instance B7}\label{graphDualB7}  
\end{figure}

\paragraph{\textbf{MIP with outage relaxations}}
We analyze now the compromise between  computation time and dual bound quality with the parametric relaxations of section \ref{sec::simpMIP}.
MIP computations are almost instantaneous for $k=1$ or $k=0$ and all the different relaxation accelerate significantly the MIP convergence compared to the exact formulation of section \ref{sec::v0MIP}.
Figure \ref{graphDualB7} illustrates the two situations which are met. 
For relatively easy instances (i.e. all instances except B7,B8,B9), the dual bounds with $k_0=5$ are very close to the exact formulation, 
with a little degradation to the target objective function (measured to convergence to optimality of the MIP relaxations).
For instances  B7, B8 and B9, a similar situation to \cite{dupin2018parallel} occurs:
although the relaxations $v_3(k_0)$ with $k_0=5$ converge to a worst lower bounds than $v_0$, the induced acceleration
of the B\&B algorithm is significant and leads to better dual solutions in limited time.
This is even true truncating the B\&B algorithm to one node, having lighter MIPs with $v_3(k_0)$ allows to generate more cuts, which is
decisive for the dual bound quality.
To explain such facts, the good primal solutions have in general $4$ or $5$ outages, a $6$th outage
induces in general over costs as T2 production is cheaper than T1 production.
It explains that $v_3(k_0)$ with $k_0=4$ or $k_0=5$ are good approximations, the major approximation with
$k_0=5$ is due to the approximation of the last penalization costs.
This possibility to have a $6$th outage induces more binaries and continuous variables, it is highly penalizing for MIP solving capacities.
This illustrates the concept of dual heuristics: relaxations are parametrized to lead to few degradation quality of the dual bounds, but
improving significantly the solving capacities.

\begin{table}[ht]
\centering
\caption{Duals bounds for the dataset A, comparison of former best dual bounds of the literature with the dual bounds $v_3(k_0)$ with $k_0>3$ and $v_0$ 
without aggregation of production time steps and decomposition scenario by scenario,
allowing at most 1h MIP computations for each scenario computation.
For $v_0$, bounds using only the LP relaxation and the MIP dual bounds without any branching are also compared.
}\label{resultDataA}
\begin{tabular}{|l|l|l|l l |l l l|l l|}
\hline
&Primal& best &$k_0=4$&$k_0=5$&$v_0$&$v_0$&$v_0$\\
&Best&\cite{Bra13}&MIP&MIP&PL&+cuts&MIP\\

 \hline
A1&169474,5M&2,31\%&0,54\%&0,08\%&0,29\%&0,09\%&0,04\%\\
A2&145956,7M&4,09\%&1,42\%&0,33\%&0,64\%&0,35\%&0,25\%\\
A3&154277,2M&3,80\%&1,85\%&0,71\%&1,06\%&0,48\%&0,37\%\\
A4&111494M&8,22\%&3,36\%&2,38\%&2,82\%&1,93\%&1,52\%\\
A5&124543,9M&9,70\%&4,25\%&3,24\%&3,46\%&2,88\%&2,55\%\\
\hline
\end{tabular}

\end{table}

\subsection{\textbf{Lower bounds for the dataset A}} For none of the instance of the dataset A, 
the hypothesis (\ref{hypCostConst}) holds to apply the time step aggregation
of section 7.2. 
Single scenario computations with disaggregated time steps are however tractable to provide dual bounds thanks to the scenario decomposition of section 7.3.
Table \ref{resultDataA} compares the dual bounds for the ROADEF challenge  computing dual bounds of $v^{sto}$
with LP relaxation, and MIP computations truncated in 1 hour,
to the bounds obtained with parametric simplified formulation of section 5 with $k^0=4$ and $k^0=5$,
improving significantly  the former best lower bounds in \cite{Bra13}.

 \begin{table}[ht]
\centering
\caption{Quick  dual bounds for datasets B and X of the EURO/ROADEF  2010 Challenge, 
with MIP computations to optimality of dual bounds $v_3(k_0)$ with $k_0<3$
with aggregation of production time steps and decomposition scenario by scenario,
compared to the ones of \cite{Bra13}}\label{quickBounds}

\begin{tabular}{|l|l|l|l|l|l|l|}
\hline
Instances&Primal&\cite{Bra13}.1&\cite{Bra13}.2&$k^0=0$&$k^0=1$&$k^0=2$\\
\hline
\hline
B6&83424,7M&56,44\%&16,58\%&19,18\%&15,12\%&11,68\%\\
B7&81109,9M&52,83\%&15,51\%&16,24\%&13,66\%&11,11\%\\
B8&81899,7M&65,31\%&23,57\%&23,79\%&20,18\%&17,03\%\\
B9&81689,5M&63,28\%&21,67\%&22,45\%&18,68\%&15,39\%\\
B10&77767M&60,92\%&18,03\%&19,58\%&16,71\%&13,43\%\\
\hline
Total B&405890,8M&59,74\%&19,08\%&20,26\%&16,87\%&13,73\%\\
\hline
X11&79007,6M&57,75\%&15,29\%&15,16\%&13,00\%&10,95\%\\
X12&77564M&52,63\%&14,19\%&14,77\%&12,40\%&10,29\%\\
X13&76288,5M&66,20\%&18,53\%&18,39\%&15,94\%&13,93\%\\
X14&76149,8M&64,67\%&17,21\%&17,61\%&15,16\%&13,08\%\\
X15&74388,4M&61,76\%&16,83\%&18,24\%&15,70\%&12,67\%\\
\hline
Total X&383398,3M&60,55\%&16,39\%&16,81\%&14,42\%&12,17\%\\
\hline
\end{tabular}


\caption{Dual bounds for datasets B and X of the EURO/ROADEF  2010 Challenge, computing of dual bounds $v_3(k_0)$ with $k_0>3$ and $v_0$ 
with aggregation of production time steps and decomposition scenario by scenario,
allowing at most 1h MIP computations for each scenario computation.
For $v_0$, bounds using only the LP relaxation and the MIP dual bounds without any branching are also compared.}\label{resultISCO}

\begin{tabular}{|l|l|l|l l |l l l|l l|}
\hline
&Primal& best &$k^0=4$&$k^0=5$&$v_0$&$v_0$&$v_0$\\
&Best&\cite{Bra13}&MIP&MIP&PL&+cuts&MIP\\
\hline
\hline
B6&83424,7M&16,58\%&6,85\%&5,36\%&7,02\%&5,72\%&\bf{5,30}\%\\
B7&81109,9M&15,51\%&7,16\%&\bf{6,23}\%&9,52\%&8,87\%&7,57\%\\
B8&81899,7M&23,57\%&14,50\%&\bf{12,37}\%&16,98\%&15,48\%&15,48\%\\
B9&81689,5M&21,67\%&12,65\%&\bf{11,59}\%&16,36\%&15,14\%&15,14\%\\
B10&77767M&18,03\%&8,33\%&\bf{7,32}\%&10,01\%&8,27\%&7,71\%\\
\hline
Total B&405890,8M&19,08\%&9,90\%&8,58\%&11,99\%&10,70\%&10,25\%\\
\hline
X11&79007,6M&15,29\%&6,32\%&5,88\%&6,60\%&5,73\%&\bf{5,45}\%\\
X12&77564M&14,19\%&7,39\%&\bf{6,65}\%&9,02\%&7,46\%&7,03\%\\
X13&76288,5M&18,53\%&9,98\%&7,90\%&8,72\%&7,92\%&\bf{7,58}\%\\
X14&76149,8M&17,21\%&9,62\%&\bf{7,82}\%&9,67\%&8,76\%&8,34\%\\
X15&74388,4M&16,83\%&9,16\%&\bf{8,11}\%&10,90\%&8,81\%&8,31\%\\
\hline
Total X&383398,3M&16,39\%&8,47\%&7,26\%&8,99\%&7,72\%&7,32\%\\
\hline
\end{tabular}

\end{table}

\subsection{\textbf{Lower bounds for the datasets B and X with scenario decomposition}}
The hypotheses (\ref{hypCostConst}) are verified for datasets B and X.
Hence, dual bounds can be computed with time step aggregation and scenario decomposition thanks to Proposition \ref{dualboundAgregTimeStep}.
The analyzes of section 8.2 ensures that computations of duals bounds are tractable for instances from datasets B and X
with LP computations or the dual bounds of truncated B\&B searches
for one scenario and time step aggregation of section 7.2 and using MIP formulations of sections 4 and 6.
Table \ref{resultISCO} compares the dual bounds for the ROADEF challenge  computing dual bounds of $v_s^{det}$
with LP relaxation, and MIP computations truncated in 1 hour,
to the bounds obtained with parametric simplified formulation of section 6 with $k^0=4$ and $k^0=5$,
to the former best bounds of the state of the art in \cite{Bra13}.
Quick dual bounds computations are also possible with $k^0\leqslant 2$ using formulations $v_3(k_0)$, as reported in Table \ref{quickBounds}.


At this stage, our dual bounds outclass significantly the former best dual bounds of the literature.
We note that the best bounds of \cite{Bra13} were slightly better than the first parametric dual bounds with $k^0=0$,
parameter $k^0=1$ already improved the bounds of \cite{Bra13}.
 Comparing the dual bounds got computing lower bounds of $v_s^{det}$ for all $s\in\SS$ emphasizes the impact of the cuts and branching of Cplex already shown with Table 1.
With $1h$ limit for MIP computations, parametric formulations with $k^0=5$ often gives  better dual bounds than the exact MIP formulation for the most difficult instances, extending the situations
met and illustrated with Figure \ref{graphDualB7}.

\subsection{\textbf{New dual bounds for the challenge EURO/ROADEF 2010}}

%


Adding  partial scenario decomposition for datasets B and X, last results were slightly improved, as reported in the 
Table \ref{bestDualROADEF}.
These last results are our best computed results.
We note that having several scenarios induces more difficult problems, especially to compute the LP relaxations.
In the search of the good trade-off between acceleration of computation
and quality of dual bounds, it was more often efficient than before to use the parametric relaxation $v_3(k_0)$ with $k_0=5$
instead of $v_0$ computations.

\begin{table}[ht]
\centering

\caption{Best dual bounds for EURO/ROADEF 2010  Challenge, comparison to the former ones of \cite{Bra13}.
$k=\infty$ denotes that the exact formulation of Section 4 was used for the best results, otherwise
$k=5$ denotes that relaxation $v_3(5)$ gave the best results.
Gaps are calculated to the best primal bounds, given mainly by \cite{gardi}
and reported in \cite{Gav13} with new best solutions for A1,A2 and A3.}\label{bestDualROADEF}

\begin{tabular}{|l|c| c c| c c |c c |}
\hline
Instances&Primal&\cite{Bra13} dual&Gap&Our best&Gap&k&nbScenar\\
&Best&Best&&duals&&&\\
\hline
\hline
A1&169474,5M&165560M&2,31\%&169403M&0,04\%&$\infty$&1\\
A2&145956,7M&139991M&4,09\%&145593M&0,25\%&$\infty$&1\\
A3&154277,2M&148454M&3,80\%&153704M&0,37\%&$\infty$&1\\
A4&111494M&102326M&8,22\%&109801M&1,52\%&$\infty$&1\\
A5&124543,9M&112467M&9,70\%&121366M&2,55\%&$\infty$&1\\
\hline
Total A&705785,5M&668798&5,24\%&699867M&0,84\%&&\\
\hline

\hline
B6&83424,7M&69592M&16,58\%&79205M&5,06\%&5&3\\
B7&81109,9M&68528M&15,51\%&76356M&5,86\%&5&3\\
B8&81899,7M&62594M&23,57\%&72023M&12,06\%&5&3\\
B9&81689,5M&63991M&21,67\%&72437M&11,33\%&5&3\\
B10&77767M&63747M&18,03\%&72378M&6,93\%&5&5\\
\hline
Total B&405890,8M&328452M&19,08\%&372399M&8,25\%&&\\
\hline
X11&79007,6M&66931M&15,29\%&74715M&5,43\%&$\infty$&5\\
X12&77564M&66558M&14,19\%&72768M&6,18\%&5&5\\
X13&76288,5M&62155M&18,53\%&70840M&7,14\%&$\infty$&3\\
X14&76149,8M&63045M&17,21\%&70373M&7,59\%&$\infty$&3\\
X15&74388,4M&61866M&16,83\%&68700M&7,65\%&5&5\\
\hline
Total X&383398,3M&320555M&16,39\%&357396M&6,78\%&&\\
\hline
\end{tabular}
\end{table}

\section{Conclusions and perspectives}

\paragraph{\textbf{Conclusions}}
New dual bounds for the 2010 EURO/ROADEF Challenge  are provided in this paper. 
The relaxation of CT6 and CT12 allows to have a MIP formulation declaring binaries only for the decisions of outage dates.
To deal with smaller problems, we proved that  dual bounds can be computed respectively aggregating production time steps and restricting to deterministic computations.
A  parametric aggregation of outages leads to a parametric family of dual bounds.
This leads to tractable computations,  outclassing significantly the former best dual bounds of the literature.
The  methodology  has similarities with meta-heuristics:
 better solutions are obtained degrading slightly the objective function, but having quicker computations and better MIP convergences.
In order to analyze the trade-off between cost degradation and acceleration of MIP computations, the knowledge of the problem is crucial.


Intermediate results have also their importance for a better understanding of the structural difficulties for the EURO/ROADEF 2010  Challenge.
The solving limits justify the simplified hypotheses and aggregation used in \cite{lusby}.
An original formulation of the stretch constraints CT6 is furnished, but it is not prominent for the quality of the
LP relaxation. 
The quality of dual bounds justifies  the quality of the primal solutions of \cite{gardi}.
Also the difficulty induced by the real world size instances for exact based approaches tend to
justify aggressive local search approaches for an industrial implementation to schedule maintenances of nuclear power plants.

\paragraph{\textbf{Perspectives}} This work offers new perspectives to improve the dual bounds.
First, using more powerful computers and more computation time  will surely improve the reported results.
With the continuous progress of computing capacities and B\&B implementations, 
 our algorithm  will have continuously better results.
Other perspectives are to derive primal matheuristics.
Dual bounds can be computed quickly, which can guide local search algorithms.
Lastly, we note special perspectives for the approach deployed by \cite{lusby}. 
The variable definition with step variables is useful for their branching phase.
Similarly with \cite{ucp}, the LP relaxation is identical with their variables, the branching phase
would use a polyhedral isomorphism to use step variables only for the branching phase.
Furthermore, the MIP relaxations of section \ref{sec::simpMIP} would be useful
to accelerate Bender's decomposition phases, to compute quicker interesting slightly underestimated cuts,
and also to stabilize the Bender's decomposition.

\subsection*{\textbf{Acknowledgements}}
The results of this paper were obtained mainly in the PhD Thesis \cite{dupin2015modelisation} financed by the 
 French Defense Procurement Agency of the French Ministry of Defense (DGA).
  Some results of this paper were presented in \cite{dupin2016dual} in the conference Matheuristics 2016, 
  the  Sixth International Workshop on Model-based Metaheuristics in Brussel
 organized by the Universit\'e Libre de Bruxelles (ULB), Belgium.

\section*{Appendix A: Proof of Propositions \ref{tightenTW} and \ref{tightenContVar}}

In the section 7.1, two propositions were mentioned reduce the number of variables in the MIP computations.
This appendix provides their proof. 



\vskip 0.3cm

\noindent{\textbf{Proof of Proposition \ref{tightenTW}}}: 
$\mathbf{Lmin}_{i,k} = \left\lceil \frac {\mathbf{Rmin}_{i,k} - \mathbf{Amax}_{i,k}} {\mathbf{D}^w  P_i}\right\rceil$
is a first lower bound for the production cycle $(i,k)$:
the minimum fuel level after refueling is at least the minimal refueling $\mathbf{Rmin}_{i,k} $ with the positivity of fuel stocks
with CT11.
This minimal length of the cycle consider a minimum refueling and a maximal fuel consumption.

Let $\widetilde{To}_{i,k}$ and $\widetilde{Ta}_{i,k}$ the minimal and maximal weeks to process  outage $(i,k)$.
Denoting $W_{i,k} = \sum_w (1-d_{i,k,w})$ the week when outage $(i,k)$ begins, we have relations
$W_{i,k+1} \geqslant W_{i,k} + \mathbf{Da}_{i,k} + \mathbf{Lmin}_{i,k}$.
$W_{i,k+1} \geqslant \widetilde{To}_{i,k} + \mathbf{Da}_{i,k} + \mathbf{Lmin}_{i,k}$.
This minoration is valid for all feasible solution, taking the lower bound of the LHS induces:
$\widetilde{To}_{i,k+1} \geqslant \widetilde{To}_{i,k} + \mathbf{Da}_{i,k} + \mathbf{Lmin}_{i,k}$.
(\ref{tightenTWo}) is thus a valid preprocessing to tighten values of $\widetilde{To}_{i,k}$.

From relations $W_{i,k+1} \geqslant W_{i,k} + \mathbf{Da}_{i,k} + \mathbf{Lmin}_{i,k}$,
we have also $\widetilde{Ta}_{i,k+1} \geqslant W_{i,k} + \mathbf{Da}_{i,k} + \mathbf{Lmin}_{i,k}$
and then $\widetilde{Ta}_{i,k+1} \geqslant \widetilde{To}_{i,k} + \mathbf{Da}_{i,k} + \mathbf{Lmin}_{i,k}$.
(\ref{tightenTWa}) is thus a valid preprocessing to tighten values of $\widetilde{Ta}_{i,k}$.

\vskip 0.3cm

\noindent{\textbf{Proof of Proposition \ref{tightenContVar}}}:
We notice that equations (\ref{PANconso}) and (\ref{PANpertes}) induce a linear system of equalities
considering $x_{i,k,s}^{init},x_{i,k,s}^{fin}$ as variables.
A recursion formula can be deduced for $x_{i,k,s}^{init}$ where the initial condition is $x_{i,0,s}^{init} = \mathbf{Xi}_{i}$:
$$x_{i,k,s}^{init} - \mathbf{Bo}_{i,k} = r_{i,k} + \frac{\mathbf{Q}_{i,k} -1}{\mathbf{Q}_{i,k}} (x_{i,k-1,s}^{init} - \sum_w \mathbf{D}^w \:  p_{i,k-1,s,w} - \mathbf{Bo}_{i,k-1})$$
This is a induction formula $x_{i,k,s}^{init} = a_{i,k-1,s} x_{i,k-1,s}^{init} + b_{i,k-1,s}$ 
with $a_{i,k,s} = \frac{\mathbf{Q}_{i,k+1} -1}{\mathbf{Q}_{i,k+1}}$ and
$b_{i,k,s} = r_{i,k+1} - \frac{\mathbf{Q}_{i,k+1} -1}{\mathbf{Q}_{i,k+1}} \left(\sum_w \mathbf{D}^w \:   p_{i,k,s,t} - \mathbf{Bo}_{i,k} \right)$.
Lemma 1 allows also to compute  $x_{i,k,s}^{init}$ as linear expressions of the variables $r_{i,k'}$ and  $p_{i,k',s,t}$ for $k'<k$.
$x_{i,k,s}^{fin}$ are also  linear expressions of the variables 
$r_{i,k'}$ and  $p_{i,k',s,t}$ for $k'<k$. reporting last inequality in (\ref{PANconso}).

\begin{lemma}
   Let $(u_n)_{n \in \NN}$ a sequence defined by induction with $u_{n+1} = a_n u_n + b_n$, 
 with $(a_n)_{n \in \NN}$, $(b_n)_{n \in \NN}$ real numbers. We have:
  \begin{small}
\begin{equation}
u_n = \left( \prod_{l=0}^{n-1} a_l \right) u_o + \sum_{l=0}^{n-1}  \left( \prod_{m=l+1}^{n-1} a_m  \right) b_l\label{lemme}\end{equation} 
\end{small}
\end{lemma}

\noindent{\textbf{Proof of Lemma 1}:  We prove the result by induction on  $n\in \NN$.
For $n=0$ or $n=1$, the initialization is trivial, with null sums or products.
Let us suppose $n\geqslant 1$ and (\ref{lemme}) true for $n-1$.
$u_{n+1} = a_n u_n + b_n$, using induction hypothesis:
}
\small{
\begin{eqnarray*}
 u_{n+1} & = &a_n \left(\left( \prod_{l=0}^{n-1} a_l \right) u_o + \sum_{l=0}^{n-1}  \left( \prod_{m=l+1}^{n-1} a_m  \right) b_l \right) + b_n\\ 
 &= &  a_n \left( \prod_{l=0}^{n-1} a_l \right) u_o +   \sum_{l=0}^{n-1} a_n  \left( \prod_{m=l+1}^{n-1} a_m  \right) b_l  + b_n\\
& = &  \left( \prod_{l=0}^{n} a_l \right) u_o +   \sum_{l=0}^{n-1} \left( \prod_{m=l+1}^{n} a_m  \right) b_l  + b_n \left( \prod_{m=n+1}^{n} a_m  \right)\\
& = & \left( \prod_{l=0}^{n} a_l \right) u_o +   \sum_{l=0}^{n} \left( \prod_{m=l+1}^{n} a_m  \right) b_l.\phantom{2} \square
 \end{eqnarray*}
}

\end{document}